\documentclass{amsart}
\usepackage{times}
\def\mbf{\mathbf}
\newtheorem{theorem}{Theorem}
\newtheorem{proposition}{Proposition}
\newtheorem{lemma}{Lemma}

\theoremstyle{definition}

\theoremstyle{remark}

\numberwithin{equation}{section}

\def\sumprime_#1{\setbox0=\hbox{$\scriptstyle{#1}$}
\setbox2=\hbox{$\displaystyle{\sum}$}
\setbox4=\hbox{${}'\mathsurround=0pt$}
\dimen0=.5\wd0 \advance\dimen0 by-.5\wd2
\ifdim\dimen0>0pt
\ifdim\dimen0>\wd4 \kern\wd4 \else\kern\dimen0\fi\fi
\mathop{{\sum}'}_{\kern-\wd4 #1}}
\newcommand{\gs}{\mathfrak{S}}
\def\ndiv{\not \hskip .03in \mid}
\begin{document}
\title{Small Gaps Between Primes I}
\author{D. A. Goldston}
\address{Department of Mathematics, San Jose
State University, San Jose, CA 95192, USA}
\email{goldston@math.sjsu.edu}
\thanks{Goldston was supported by NSF; Y{\i}ld{\i}r{\i}m was supported
by T\"{U}B\.{I}TAK}
\author{\mbox{~} C. Y. Yildirim}
\address{ Department of Mathematics, Bo\~{g}azi\c{c}i University,
Bebek, Istanbul, 34342 Turkey \& \newline \hspace*{.3cm}
Feza G\"{u}rsey Enstit\"{u}s\"{u}, \c{C}engelk\"{o}y,
Istanbul, P.K. 6, 81220 Turkey}
\email{yalciny@boun.edu.tr}
\subjclass{Primary 11N05 ; Secondary 11P32}

\date{\today}

\keywords{prime number}

\begin{abstract} We use  short divisor sums to approximate prime tuples and moments for primes in short intervals. By connecting these results to classical moment problems we are able to prove that, for any $\eta>0$,  a positive proportion of consecutive primes are within $\frac{1}{4} + \eta$ times the average spacing between primes.
 \end{abstract}

\maketitle

\section{Introduction}
Finding mathematical proofs for easily observed properties of the distribution of prime numbers is a difficult and often humbling task, at least for the authors of this paper.  The twin prime conjecture is a famous example of this, but we are concerned here with the much more modest problem of proving that there are arbitrarily large primes that are \lq \lq unusually close "  together. Statistically this means we seek consecutive primes whose distance apart is substantially less than the average distance between consecutive primes. Letting  $p_n$ denote the $n^{\mathrm{th}}$ prime, then by the prime number theorem the average gaps size $p_{n+1}-p_n$ between consecutive primes is $\log p_n$. Thus we define 
\begin{equation}  \Delta  =  \liminf_{n\to \infty} \left( \frac{p_{n+1}-p_n}{\log p_n}\right),\label{1.1} 
\end{equation}
so that $\Delta$ is the smallest number for which there will be infinitely many  gaps between consecutive primes of size less than $\Delta +\epsilon$ times the average size. It is empirically evident that
\begin{equation}  \Delta =0, \label{1.2} \end{equation}
but this has never been proved.  Up to now three different unconditional methods have been invented which provide non-trivial estimates for $\Delta$.

\bigskip
\noindent \textbf{1. The Hardy-Littlewood and Bombieri-Davenport Method.} In the mid-1920's Hardy and Littlewood \footnote{In the unpublished paper Partitio Numerorum VII they proved,  assuming the Generalized Riemann Hypothesis,  that $\Delta\le \frac {2}{3}$.  In 1940 Rankin \cite{Ra} refined Hardy and Littlewood's method to show that 
$\Delta \le \frac {1 + 4\Theta}{5}$, where $\Theta$ is the supremum of the real parts of all the zeros of all Dirichlet $L$-functions. In particular assuming the Generalized Riemann Hypothesis ($\Theta = \frac{1}{2}$) this gives $\Delta \le \frac{3}{5}$.} used the circle method to obtain a conditional result which in 1965 Bombieri and Davenport \cite{BD} both improved and made unconditional. This approach can be interpreted as a second moment method using a truncated divisor sum as an approximation of $\Lambda(n)$, the von Mangoldt function (see introduction in \cite{GYI}). The method proves 
\begin{equation}  \Delta \le \frac{1}{2}. \label{1.3} \end{equation}

\bigskip
\noindent \textbf{2. The Erd\"os Method.} By the prime number theorem we have $\Delta \le 1$.  Erd\"os \cite{E} in 1940 was the first to prove unconditionally that $\Delta <1$. He used the sieve upper bound for primes differing by an even number $k$ 
\begin{equation} \sum_{n\le N}\Lambda(n)\Lambda(n+k) \le (\mathcal{B}+\epsilon)\gs(k) N \label{1.4}\end{equation} 
where $\gs(k)$ is the singular series and $\mathcal{B}$ is a constant.  By this bound there can not be too many pairs of primes with the same difference, and therefore the distribution function for prime gaps must spread out from the average. This method gives the result
\begin{equation} \Delta \le 1 - \frac{1}{2\mathcal{B}}. \label{1.5} \end{equation} 
The value $\mathcal{B}=4$ of Bombieri and Davenport \cite{BD} (see also \cite{Gr} or \cite{HR}) or $\mathcal{B}=3.5$  of Bombieri, Friedlander, and Iwaniec \cite{BFI}\footnote{This value only holds for $k$ not too large as a function of $N$ in \eqref{1.4}, but this is acceptable for \eqref{1.5}.}, or even slightly smaller values  may be used here.

\bigskip
\noindent \textbf{3. The Maier Method.} In 1988 Maier \cite{Ma} found certain intervals (rather sparsely distributed) where there are $e^\gamma$ more primes than the expected number, and therefore within these intervals the average spacing is reduced by a factor of $e^{-\gamma}$. Hence
\begin{equation} \Delta \le e^{-\gamma} = 0.56145\ldots  \ .\label{1.6} \end{equation}
In contrast to the first two methods, this method does not produce a positive proportion of small prime gaps. 

\bigskip
These methods may be combined to obtain improved results. Huxley \cite{Hu1,Hu2} combined the first two methods making use of a weighted version of the first method to find 
\begin{equation}  \Delta \le 0.44254\ldots \quad (\mathrm{using} \ \mathcal{B}=4), \quad \quad  
\Delta \le 0.43494\ldots \quad (\mathrm{using }\ \mathcal{B}=3.5), \label{1.7} \end{equation}
and Maier combined his method with Huxley's method with $\mathcal{B} =4$ to obtain
\begin{equation} \Delta \le ( 0.44254\ldots )e^{-\gamma} = .24846 \ldots \ .\label{1.8} \end{equation}
This last result is the best result known up to now, and as we have seen uses all three of the previously known methods.

\bigskip
 For several years we have been developing tools for dealing with higher correlations of short divisor sums which approximate primes. Our first results appeared in \cite{GYI}, and, with considerable help from other mathematicians, we have greatly simplified and improved on these results in \cite{GYIII}. In the former paper we had an application to small gaps between primes based on approximating a third moment. In particular we recovered the result \eqref{1.3}. The method is based on the same approximation that underlies the method of Bombieri-Davenport, but it detects small prime gaps in a different way. In this paper we extend that argument to all moments and obtain the limit of this method. 

Let $\pi(N)$ denote the number of primes less than or equal to $N$. 
\begin{theorem}  Let $r$ be any positive integer. For any fixed $\lambda > ( \sqrt{r}-\frac{1}{2})^2 $ and $N$ sufficiently large, we have
\begin{equation} \sum_{\substack{ p_n \le N\\ p_{n+r}-p_n \le \lambda \log p_n}}1\  \gg_r \pi(N). \label{1.9}\end{equation}
In particular, for any fixed $\eta >0$ and all sufficiently large $N>N_0(\eta)$, a positive proportion of gaps $p_{n+1}-p_n$ for $p_n\le N$ are less than $(\frac{1}{4} + \eta)\log N$, and 
\begin{equation}  \Delta \le \frac{1}{4}. \label{1.10} \end{equation}

\end{theorem}

Our results depend on the level of distribution of primes in arithmetic progressions, and Theorem 1 makes use of the Bombieri-Vinogradov theorem. If for primes up to $N$ the level of distribution in arithmetic progressions is assumed to be $N^{\vartheta-\epsilon}$ for any $\epsilon>0$, then Theorem 1 holds with $\lambda > (\sqrt{r} - \sqrt{\frac{\vartheta}{2}}\ )^2$. Hence, assuming the Elliott-Halberstam conjecture that $\vartheta=1$ holds, we obtain the improved result that 
\begin{equation} \liminf_{n\to \infty} \left( \frac{p_{n+r}-p_n}{\log p_n}\right) \le  \left( \sqrt{r} - \frac{1}{\sqrt{2}} \right)^2, \label{1.11}
\end{equation}
and in particular
\begin{equation} \Delta \le \left(\frac{3}{2}-\sqrt{2}\right) = 0.085786\ldots \ < \frac{1}{11}. \label{1.12}
\end{equation}

There are several improvements that can be made in our results. First,  we can incorporate Maier's method into our method. This is a straightforward adaptation of the argument Maier used to combine his method with Huxley's result, although the result is complicated by the need to prove our propositions in the next section when they are summed over arithmetic progressions.  Second, and more significantly, we have found in joint work with J. Pintz better approximations for prime tuples than those used in this paper, and these lead to significantly stronger results. These results will appear in future papers. 

This paper is organized as follows.  In section 2 we present our method and state the two main propositions needed in the proof. In section 3 we prove some lemmas which are used in the later sections. In sections 4 and 5 we prove the propositions. In section 6 we examine an optimization problem related to the Poisson distribution which is used in the proof of Theorem 1, and finally in  section 7 we prove  Theorem 1. 

\emph{ Notation.}  We will take $\epsilon$ to be any sufficiently small positive number whose value can be changed from equation to equation, and similarly $C$, $c$, and $c'$ will denote  small fixed positive constants whose value may change from equation to equation. We will let $A$ denote a large positive constant which may be taken as large as we wish, but is fixed throughout the paper. For a finite set $\mathcal{A}$ we let $ |\mathcal{A}|$ denote the number of elements in $\mathcal{A}$. We will sometimes write $\mathcal{A}=\mathcal{A}_k$ if $|\mathcal{A}|=k$. For a vector $\mbf{H}$ we denote the number of components by $|\mbf{H}|$. A dash in a summation sign $\sum '$ indicates all the summation variables are relatively prime with each other, and further any sum without a lower bound on the summation variables will have the variables start with the value 1. Empty sums will have the value zero, and empty products will have the value 1.  We will make use of the Iverson notation that, for a statement $P$, $[P]$ is 1 if $P$ is true and is 0 if $P$ is false. 

\emph{ Acknowledgment.} We are indebted to Andrew Granville and Kannan Soundararajan who have greatly clarified and simplified our method. We have made extensive use of their results here. We have also benefited  from ideas of Enrico Bombieri, Brian Conrey, Percy Deift, David Farmer, John Friedlander, Roger Heath-Brown, Hugh Montgomery, Michael Rubinstein, and Peter Sarnak. We have used work of Dashiell Fryer, an undergraduate student at San Jose State University in the MARC program, who worked on properties of Laguerre polynomials needed in our method. The first-named author also thanks the American Institute of Mathematics where much of the collaboration mentioned above took place. In a recent preprint \cite{GT} Ben Green and Terence Tao proved a landmark result on arithmetic progressions of primes. One tool they used was the current Proposition 1 from an earlier (not widely distributed) preprint of this paper. They corrected an oversight in our original proof which we have incorporated into our Lemma 2 and the proof of Proposition 1.

\section{Approximating Prime Tuples}

Our approach for finding small gaps between primes is to compute approximations of the moments for the number of primes in short intervals, and this computation uses short divisor sums to approximate prime tuples. Given a positive integer $h$, let 
\begin{equation} \mathcal{H}=\{h_1,h_2, \ldots, h_k\}, \quad \textrm{with } 0\le h_1,h_2, \cdots , h_k \le h \textrm{ distinct integers},\label{2.1}\end{equation}
and let $\nu_p(\mathcal{H})$ denote the number of distinct residue classes modulo $p$ the elements of $\mathcal{H}$ occupy. We define the singular series
\begin{equation}  \gs(\mathcal{H}) = \prod_p\left( 1-\frac{1}{p}\right)^{-k}\left(1 - \frac{\nu_p(\mathcal{H})}{p}\right) .\label{2.2}\end{equation}
If $\gs(\mathcal{H})\neq 0$ then $\mathcal{H}$ is called \emph{admissible}. Thus $\mathcal{H}$ is admissible if and only if $\nu_p(\mathcal{H})<p$ for all $p$. 

Letting $\Lambda(n)$ denote the von Mangoldt function, define
\begin{equation}  \Lambda(n;\mathcal{H}) = \Lambda(n+h_1)\Lambda(n+h_2)\cdots \Lambda(n+h_k). \label{2.3}\end{equation}
The Hardy-Littlewood prime tuple conjecture \cite{HL} states that for $\mathcal{H}$ admissible, 
\begin{equation}  \sum_{n\le N}\Lambda(n;\mathcal{H})  = N\big(\gs(\mathcal{H})+o(1)\big), \quad  \mbox{as \ $N\to \infty$.}\label{2.4}\end{equation}
(This is trivially true if $\mathcal{H}$ is not admissible.) We approximate $\Lambda(n)$ as in our earlier work by using the truncated divisor sum
\begin{equation} \Lambda_R(n) = \sum_{\substack{d|n\\ d\le R}}\mu(d)\log \frac{R}{d}, \label{2.5} \end{equation}
and then approximate $\Lambda(n;\mathcal{H})$ by  
\begin{equation}  \Lambda_R(n;\mathcal{H})= \Lambda_R(n+h_1)\Lambda_R(n+h_2)\cdots \Lambda_R(n+h_k). \label{2.6}\end{equation}
For convenience we also define $\Lambda_R(n;\mathcal{H})=1$ if $\mathcal{H}= \emptyset$. 
Our method is founded on the following two propositions which allow us to obtain information about primes. Suppose $\mathcal{H}_1$ and $\mathcal{H}_2$ are both sets of distinct positive integers that are $\le h$, with $|\mathcal{H}_1|=k_1$ and $|\mathcal{H}_2|=k_2$, and let $k=k_1+k_2$. We always assume $k\ge 1$. 
\begin{proposition}   Let $\mathcal{H}=\mathcal{H}_1 \cup \mathcal{H}_2$, and $r=|\mathcal{H}_1\cap \mathcal{H}_2|$. If $R=o(N^{\frac{1}{k}}(\log R)^{1-\frac{|\mathcal{H}|}{k}}))$ and $h\le R^A$ for any large constant $A>0$, then we have for  $R,N\to \infty$, 
\begin{equation} \sum_{n\le N} \Lambda_R(n;\mathcal{H}_1)\Lambda_R(n;\mathcal{H}_2) = N\big(\gs(\mathcal{H})+o_k(1)\big)(\log R)^{r} .\label{2.7}\end{equation}
\end{proposition} 
\begin{proposition} Let $\mathcal{H}=\mathcal{H}_1 \cup \mathcal{H}_2$, , $r=|\mathcal{H}_1\cap \mathcal{H}_2|$, and $1\le h_0\le h$. Let $\mathcal{H}_0 = \mathcal{H}\, \cup \, \{h_0\}$, and $r_0=r$ if $h_0\not \in \mathcal{H}$ and $r_0=r+1$ if $h_0 \in \mathcal{H}$. If $R\ll_k N^{\frac{1}{2k}}(\log N)^{-B(k)}$ for a sufficiently large positive constant $B(k)$,  and $h\le R^{\frac{1}{2k}}$, then we have for $R,N\to \infty$, 
\begin{equation} \sum_{n\le N} \Lambda_R(n;\mathcal{H}_1)\Lambda_R(n;\mathcal{H}_2)\Lambda(n+h_0) =  N\big(\gs(\mathcal{H}_0)+o_k(1)\big)(\log R)^{r_0}. \label{2.8}\end{equation}
Assuming the Elliott-Halberstam conjecture, then equation \eqref{2.8} holds for  $R\ll_k N^{\frac{1}{k} - \epsilon}$ with any $\epsilon>0$.
\end{proposition}
The restriction on the size of $R$ in Proposition 2 can be improved in the situation when $h_0 \in \mathcal{H}_1\cup \mathcal{H}_2$. If we let $k^*= k-|\mathcal{H}_1\cap \{h_0\}|-|\mathcal{H}_2\cap \{h_0\}|$, then we see that the reduction of cases at the start of the proof of Proposition 2  implies that Proposition 2 holds in the range $R\ll_k N^{\frac{1}{2k^*}}(\log N)^{-B(k)}$ except in the trivial case when $k=2$ and $k^*=0$ where the result holds for $R\le N$. 
In the case of the Elliott-Halberstam conjecture $k$ can also be replaced by $k^*$ in the bound for $R$. 

We actually prove both Propositions with the error term $o_k(1)$ replaced by a series of lower order terms, which however are not needed in any of our applications. 

If we take $\mathcal{H}_2= \emptyset$ in Proposition 1 we have, for $R=o( N^{\frac{1}{k}})$ and $h\le R^A$ for any given constant $A>0$, that for  $R,N\to \infty$, 
\begin{equation} \sum_{n\le N} \Lambda_R(n;\mathcal{H}) = N\big(\gs(\mathcal{H})+o_k(1)\big),\label{2.9}\end{equation}
in agreement with the Hardy-Littlewood prime tuple conjecture \eqref{2.4}.

In applying these propositions it is critical to have some form of positivity in the argument. For example, in the special case when $\mathcal{H}_2=\emptyset$ Proposition 2 takes the form, for $R\le N^{\frac{1}{2k}}(\log N)^{-B(k)}$,
\begin{equation} \sum_{n\le N} \Lambda_R(n;\mathcal{H})\Lambda(n+h_0) =  \left\{ \begin{array}{ll} N\big(\gs(\mathcal{H})+o_k(1)\big)\log R,  & \mbox{if \ $h_0 \in  \mathcal{H}$, }\\ N\big(\gs(\mathcal{H}_0)+o_k(1)\big) , & \mbox{if \ $h_0 \not \in  \mathcal{H},$} \end{array} \right. \label{2.10}\end{equation}
which would seem to exhibit that our approximation detects primes. However since $\Lambda_R(n;\mathcal{H})$ is not non-negative, it is impossible to conclude anything about primes from \eqref{2.10} alone. 

On the other hand, consider instead the special case of Proposition 2 where $\mathcal{H}_1=\mathcal{H}_2 = \mathcal{H}$ which gives, on taking $|\mathcal{H}|=k$, for $R\le N^{\frac{1}{4k}}(\log N)^{-B(k)}$,
\begin{equation} \sum_{n\le N} \Lambda_R(n;\mathcal{H})^2\Lambda(n+h_0) =  \left\{ \begin{array}{ll} N\big(\gs(\mathcal{H})+o_k(1)\big)(\log R)^{k+1} ,  & \mbox{ if \ $h_0 \in  \mathcal{H}$, }\\ N\big(\gs(\mathcal{H}_0)+o_k(1)\big)(\log R)^k , & \mbox{if \ $h_0 \not \in  \mathcal{H}$ }. \end{array} \right. \label{2.11}\end{equation}
The restriction on the size of $R$ here makes it impossible to conclude from \eqref{2.11} that any given tuple $\mathcal{H}$ will contain two or more primes (see \cite{GO1}), but Granville and Soundararajan found a simple argument which uses the non-negativity of $\Lambda_R(n;\mathcal{H})^2$ to prove from \eqref{2.11} that 
\begin{equation}  \Delta \le \frac{3}{4} .\label{2.12}\end{equation}
To prove their result, we need a formula of  Gallagher  that as $h \to \infty$,\footnote{Granville (unpublished) and Montgomery-Soundararajan \cite{MS} have recently proved more precise results, but these are not needed here.}  
\begin{equation} \sum_{\substack{1\le h_1,h_2,\ldots ,h_k\le h\\ \text{distinct}}}\gs(\mathcal{H})=h^k +O_k(h^{k-\frac{1}{2}+\epsilon}).\label{2.13}\end{equation}
We fix $k \ge 1$; the argument works equally well for any $k$, and we can take $k=1$ if we wish. Suppose now that 
\[ R= N^{\frac{1}{4k}}(\log N)^{-B(k)},\qquad h\ll \log N.\]
 By differencing, equation \eqref{2.11} continues to hold when the sum on the left-hand side is over $N< n\le 2N$, and therefore we have
\[ \begin{split}\sum_{n= N+1}^{2N}\Big(\sum_{1\le h_0\le h}&\Lambda(n+h_0)\Big) \Lambda_R(n;\mathcal{H})^2 \\ & \sim  kN\gs(\mathcal{H})(\log R)^{k+1}   +\sum_{\substack{1\le h_0\le h \\ h_0\neq h_i,\ 1\le i\le k}}N\gs(\mathcal{H}_0)(\log R)^{k} . \end{split}\]
Also by Proposition 1 
\[ \sum_{n=N+1}^{2N} \Lambda_R(n;\mathcal{H})^2 \sim  N\gs(\mathcal{H})(\log R)^k , \]
and therefore we have on summing over all distinct tuples $1\le h_1,h_2,\ldots ,h_k\le h$ and applying \eqref{2.13} that, for $\rho$ a fixed number and $h,N\to \infty$, 
\[ \begin{split}\sum_{n=N+1}^{2N}\Big(\sum_{1\le h_0\le h}\Lambda(n+h_0)- &\rho \log N \Big) \Big(\sum_{\substack{1\le h_1,h_2,\ldots ,h_k\le h\\ \text{distinct}}}\Lambda_R(n;\mathcal{H})^2 \Big) \\&\sim Nh^k(\log R)^k \big( k\log R +h - \rho\log N\big)\\
&\sim  N h^k(\log R)^k\big(h - (\rho - \frac{1}{4})\log N \big). \end{split}\]
Since $\Lambda_R(n)\le d(n)\log R \ll n^\epsilon$, we see that the contribution in the sum above from terms where $n+h_0$ is a prime power is $\ll N^{1/2+\epsilon}$ which is negligible, and therefore we may restrict the sum over $h_0$ to terms where $n+h_0$ is prime.  
The right-hand side above is positive if $h> (\rho - \frac{1}{4})\log N $, which implies with this restriction on $h$ that there is a value of $n$, $N<n\le 2N$, such that
\[ \sum_{\substack{1\le h_0\le h\\ n+h_0\ \text{prime}}}\log(n+h_0)> \rho \log N  .\]
If $\rho >1$ this implies that for $N$ sufficiently large there are at least two terms in this sum, and thus by taking $\rho \to 1^+$ we obtain \eqref{2.12}.

The proof of Theorem 1 is a refinement of the above argument, where we detect primes by the square of the linear combination of tuple approximations 
\begin{equation} a_0+ \sum_{j=1}^k a_j\Big(\sum_{\substack{1\le h_1,h_2,\ldots ,h_j\le h\\ \text{distinct}}}\Lambda_R(n;\mathcal{H}_j)\Big).\label{2.14} \end{equation}
Here the $a_j$'s are available to optimize the argument. While it is possible to use \eqref{2.14} directly, we have chosen in the proof of Theorem 1 to first approximate moments, which highlights the Poisson model which the prime numbers are thought to satisfy. This method also has the advantage of simplifying the combinatorics that occur in the problem. The moment method leads to an optimization problem which is familiar in the theory of orthogonal polynomials, the solution of which was provided to us by Enrico Bombieri and Percy Deift. The final result we obtain depends on the asymptotics of the smallest zero of a certain sequence of Laguerre polynomials; these  results are obtained by Sturm comparison type theorems and have appeared in the literature; Michael Rubinstein first pointed these out to us.

\section{Lemmas}

We recall the Riemann zeta-function  defined for $\mathrm{Re}(s)>1$ by
\begin{equation} \zeta(s) = \prod_p\left( 1 - \frac{1}{p^s}\right)^{-1}. \label{3.1}\end{equation}
The zeta function is analytic everywhere except for a simple pole with residue $1$ at $s=1$, and therefore
\begin{equation} \zeta(s) - \frac{1}{s-1}  \label{3.2}\end{equation}
 is an entire function. 
We need to use a classical zero-free region result.  By  Theorem 3.11 and (3.12.8) of \cite{T}  there exists a small positive constant $C$ such that $\zeta(\sigma +it)\neq 0$ in the region  
\begin{equation} \sigma \ge 1 -\frac{C}{ \log(|t|+2)} \label{3.3}\end{equation}
 for all $t$, and further 
\begin{equation} \zeta(\sigma +it)-\frac{1}{\sigma -1 +it} \ll \log(|t|+2), \hskip .4in
\frac{1}{ \zeta(\sigma +it)} \ll \log(|t|+2),
\label{3.4}
\end{equation}
in this region. 
 Let $(c)$ denote the contour $s = c + it$, $-\infty < t<\infty$, and let $\mathcal{ L}$ denote the contour given by 
\begin{equation} s= -\frac{C}{ \log(|t|+2)} +it.\label{3.5}\end{equation}
\begin{lemma} We have, for $R\ge 2$ and $c>0$ 
\begin{equation} \frac {1}{ 2\pi i}\mathop{\int}_{(c)\ } \frac{1}{\zeta(1+s)} \frac{R^{s}}{ {s}^2}\, ds = 1 + O( e^{-c'\sqrt{\log R}}), \label{3.6}\end{equation}
and for any fixed constant $B$ 
\begin{equation} \int_{\mathcal{L}}\big(\log(|s|+2)\big)^B\left|  \frac {R^s ds}{s^2}  \right| \ll e^{-c'\sqrt{\log R}}, \label{3.7} \end{equation}
and
\begin{equation} \mathop{\int}_{(\frac{1}{\log R})}\big(\log(|s|+2)\big)^B\left|  \frac {R^s ds}{s^2}  \right| \ll \log R. \label{3.8} \end{equation}
\end{lemma}
\bigskip
\noindent \emph{Proof.} We first prove \eqref{3.7}. The integral to be bounded is, for any $w\ge 2$,
\[ \begin{split} &\ll \int_{-\infty}^\infty R^{-\frac{C}{ \log(|t|+2)}}\frac{\big(\log(|t|+2)\big)^B}{(|t|+C)^2}\, dt \\ 
& \ll (\log w)^B \int_0^w R^{-\frac{C}{\log(t+2)}}\, dt + \int_w^\infty \frac{(\log t)^B}{ t^2}\, dt \\
 &\ll (w(\log w)^B) e^{\frac{-C\log R}{ \log w}} + \frac{(\log w)^B }{ w} ,\end{split}\]
and on choosing $\log w = \frac{1}{ 2}\sqrt{C\log R}$ this is 
\[ \ll  (C\log R)^{\frac{B}{2}}
 e^{-\frac{1}{ 2} \sqrt{C\log R}} \ll e^{-c'\sqrt{\log R}},\]
which proves \eqref{3.7}.

To prove \eqref{3.6}, we note that by the second bound in \eqref{3.4} 
 the integrand in \eqref{3.6}  vanishes as $|t| \to \infty$ in the region to the right of $\mathcal{L}$, and therefore we can move the contour $(c)$ to the left to $\mathcal{L}$, pass the simple pole at $s=0$ with residue $1$, and obtain
\[\frac {1}{ 2\pi i}\mathop{\int}_{(c)\ } \frac{1}{\zeta(1+s)} \frac{R^{s}}{ {s}^2}\, ds = 1 + \frac {1}{ 2\pi i}\mathop{\int}_{(\mathcal{L})\ } \frac{1}{\zeta(1+s)} \frac{R^{s}}{ {s}^2}\, ds .\]
On $\mathcal{L}$ we have from \eqref{3.4} that $\frac{1}{\zeta(1+s)}\ll \log(2+|t|)$, and therefore we may use the estimate \eqref{3.7} to obtain \eqref{3.6}. Finally, the left-hand side of \eqref{3.8} is
\[ \begin{split} &\ll \int_{|t|\le \frac{1}{\log R} } (\log R)^2 dt + \int_{|t|> \frac{1}{\log R} }\frac{\big(\log(|t|+2)\big)^B}{t^2}dt \\ & \ll \log R . \end{split}\]

\begin{lemma} Let $f_R(s_1,s_2)$ be analytic in the strip $-B\le \sigma_1,\sigma_2 \le b$ for some positive constants $B$ and $b$, and suppose also that for any $\epsilon>0$  $f_R(s_1,s_2)\ll e^{\epsilon \sqrt{\log R}}$ in this strip as $|t_1|,|t_2| \to \infty$. For $R\ge 2$ and $0<c_1,c_2 \le c$, let
\begin{equation} \mathcal{U}(R) = \frac {1}{ (2\pi i)^{2}}\mathop{\int }_{(c_2)\ } \!\mathop{\int}_{(c_1)\ } f_R(s_1,s_2) \frac{\zeta(1+s_1+s_2)}{\zeta(1+s_1)\zeta(1+s_{2})}\frac{R^{s_1+s_2}}{{s_1}^2{s_2}^2}ds_1\, ds_2 . \label{3.9}
\end{equation}
Then
\begin{equation} \mathcal{U}(R) = f_R(0,0)\log R + \mathcal{C}_R + O(e^{-c'\sqrt{\log R}}), \label{3.10} \end{equation}
where 
\begin{equation} \mathcal{C}_R = \frac{\partial f_R}{\partial s_2}(0,0) + \frac {1}{ 2\pi i}\mathop{\int}_{\mathcal{L}\ } f_R(s_1,-s_1) \frac{ds_1}{\zeta(1+s_1)\zeta(1-s_{1}){s_1}^4}.
\label{3.11} \end{equation}
\end{lemma}
\bigskip
\noindent \emph{Proof.} We write $f=f_R$ and $\mathcal{U}=\mathcal{U}(R)$ in what follows.
To evaluate $\mathcal{U}$, we first move $(c_1)$ to $\mathcal{L}$ passing a simple pole at $s_1=0$ and obtain
\begin{equation} \begin{split} \mathcal{U} &= \frac {1}{ 2\pi i}\mathop{\int }_{(c_{2})\ } f(0,s_2) \frac{R^{s_2}}{{s_2}^2} ds_2\\& \hskip .5in + 
\frac {1}{ (2\pi i)^{2}}\mathop{\int }_{(c_{2})\ } \!\mathop{\int}_{\mathcal{L}\ } f(s_1,s_2) \frac{\zeta(1+s_1+s_2)R^{s_1+s_2}}{\zeta(1+s_1)\zeta(1+s_{2}){s_1}^2{s_2}^2}ds_1\, ds_2 \\& = \mathcal{U}_1+\mathcal{U}_2,\end{split}\label{3.12}\end{equation}
where we have taken $c_2> \frac{C}{\log 2}$ so that we did not pass the simple pole of $\zeta(1+s_1+s_2)$ at $s_1=-s_2$ in $\mathcal{U}_2$. To evaluate $\mathcal{U}_1$ we move $(c_2)$ to $\mathcal{L}$ passing a double pole at $s_2=0$ and obtain by \eqref{3.7} of Lemma 1 and the bound for $f$ that
\begin{equation} \begin{split} \mathcal{U}_1 &= f(0,0)\log R + \frac{\partial f}{\partial s_2}(0,0) + \frac {1}{ 2\pi i}\mathop{\int }_{\mathcal{L}\ } f(0,s_2) \frac{R^{s_2}}{{s_2}^2} ds_2 \\& = f(0,0)\log R + \frac{\partial f}{\partial s_2}(0,0) + O(e^{-c'\sqrt{\log R}}) .\end{split}\label{3.13}\end{equation}
For  $\mathcal{U}_2$ we move $(c_2)$ to $\mathcal{L}$, passing simple poles at $s_2=-s_1$ and $s_2=0$ and obtain
\begin{equation} \begin{split} \mathcal{U}_2 &= 
\frac {1}{ 2\pi i}\mathop{\int}_{\mathcal{L}\ } f(s_1,-s_1) \frac{ds_1}{\zeta(1+s_1)\zeta(1-s_{1}){s_1}^4}+ \mathop{\int}_{\mathcal{L}\ } f(s_1,0) \frac{R^{s_1}}{{s_1}^2}ds_1\\& \qquad +\!\mathop{\int}_{\mathcal{L}\ }\!\mathop{\int}_{\mathcal{L}\ } f(s_1,s_2) \frac{\zeta(1+s_1+s_2)R^{s_1+s_2}}{\zeta(1+s_1)\zeta(1+s_{2}){s_1}^2{s_2}^2}ds_1\, ds_2 \\&  = \frac {1}{ 2\pi i}\mathop{\int}_{\mathcal{L}\ } f(s_1,-s_1) \frac{ds_1}{\zeta(1+s_1)\zeta(1-s_{1}){s_1}^4}+ O(e^{-c'\sqrt{\log R}}),  \end{split}\label{3.14}\end{equation}
where we estimated the last two integrals using \eqref{3.4}, Lemma 1, and the bound for $f$ stated in the lemma. This completes the proof. 

\section{Proof of Proposition 1} 
Let 
\begin{equation} \kappa = |\mathcal{H}_1 \cup \mathcal{H}_2|,\quad r=|\mathcal{H}_1 \cap \mathcal{H}_2|, \quad k = |\mathcal{H}_1| + |\mathcal{H}_2| =k_1+k_2, \label{4.1} \end{equation} and therefore $0\le r,k_1,k_2 \le \kappa$ and
\begin{equation} \kappa = k-r .\label{4.2} \end{equation}
Next, without loss of generality, we take
\begin{equation} \begin{split}& \mbf{H} = (h_1,h_2,\ldots , h_k), \\&
\mathcal{H}_1 = \{ h_1,h_2,\ldots , h_{k_1}\}, \\&
\mathcal{H}_2 = \{ h_{k_1+1},h_{k_1+2},\ldots , h_{k}\}, \\&
\mathcal{H}_1 \cap \mathcal{H}_2 = \{ h_1,h_2, \ldots , h_r\}, \\ & 
h_1=h_k, \quad h_2 =h_{k-1},\quad \cdots , h_r = h_{k-r+1}=h_{\kappa+1},\\ &
\mathcal{H} := \mathcal{H}_1 \cup \mathcal{H}_2 =\{h_1,h_2, \ldots , h_\kappa \}.  \end{split}\label{4.3} \end{equation}
Here  $r=0$ when $\mathcal{H}_1 \cap \mathcal{H}_2= \emptyset$ and the fourth and fifth lines in \eqref{4.3} may be removed.  With this notation we have 
\begin{equation} \begin{split} \mathcal{S}_k(\mathcal{H}_1,\mathcal{H}_2)  & = \sum_{n\le N} \Lambda_R(n;\mathcal{H}_1)\Lambda_R(n;\mathcal{H}_2)\\& = \sum_{n\le N}\sum_{\substack{d_1,d_2,\ldots,d_k\le R\\  d_i|n+h_i, 1\le i\le k}}  \prod_{i=1}^{k}\mu(d_i)\log \frac{R}{d_i}
\\& = \sum_{d_1,d_2,\ldots,d_k\le R} \Big( \prod_{i=1}^{k}\mu(d_i)\log \frac{R}{d_i}\Big)\sum_{\substack{n\le N\\  d_i|n+h_i, \ 1\le i\le k}} 1
.\end{split} \label{4.4}\end{equation}
Let
\begin{equation}  D_k=[d_1,d_2,\ldots , d_k], \label{4.5}\end{equation}
the least common multiple of $d_1, d_2, \ldots , d_k$.
The sum over $n$ above is zero unless $(d_i,d_j)|h_j-h_i$, $1\le i<j\le k$, in which case the sum runs through a unique residue class modulo $D_k$, and we have 
\begin{equation}  \sum_{\substack{n\le N\\ d_j |n+h_j, 1\le j\le k }}1= \frac{N}{D_k}+O(1). \label{4.6} \end{equation}
Hence 
\begin{equation}\begin{split} \mathcal{S}_k(\mathcal{H}_1,\mathcal{H}_2) &= N\sum_{\substack{d_1,d_2,\cdots d_k\le R\\(d_i,d_j)|h_j-h_i,1\le i<j\le k}} \frac{1}{D_k}\prod_{j=1}^k\mu(d_j)\log \frac{R}{d_j} +O(R^k)\\
&= N\mathcal{T}_k(\mathcal{H}_1,\mathcal{H}_2) +O(R^k).\label{4.7}\end{split}\end{equation}
We next decompose the $d_i$'s into relatively prime factors. Let $\mathcal{P}(k)$ be the set of all non-empty subsets of the set of $k$ elements $\{ 1, 2, \ldots , k\}$ (This is just the power set with the empty set removed.) For 
$\mathcal{B} \in \mathcal{P}(k)$, we let $\mathcal{P}_{\mathcal{B}}(k)$ denote the set of all members of $\mathcal{P}(k)$ for which $\mathcal{B}$ is a subset. Thus for example
if $k= 4$ then 
\[\mathcal{P}_{\{1,2\}}(4) = \{ \{1,2\}, \{1,2,3\}, \{1,2,4\}, \{1,2,3,4\} \} \]
Since the $d_i$'s are squarefree we can decompose them into the relatively prime factors 
\begin{equation}  d_i = \prod_{\nu \in \mathcal{P}_{\{i\}}(k)} a_\nu ,
\qquad 1\le i \le k ,\label{4.8}\end{equation}
where $a_\nu$ is the product of all the primes that precisely divide all the $d_i$'s for which $i \in \nu$, and none of the other $d_i$'s.  This decomposition is unique and the $2^k-1$ $a_\nu$'s are pairwise relatively prime with each other.  

We next denote by $\mathcal{D}(\mbf{H})$ the divisibility conditions
\begin{equation} (d_i,d_j)= \prod_{\nu\in \mathcal{P}_{ \{i,j\}}(k)}a_\nu \ \Big|\  h_j-h_i, \qquad 1\le i<j\le k , \label{4.9}\end{equation}
and have
\begin{equation} \mathcal{T}_k(\mathcal{H}_1,\mathcal{H}_2) = \sum_{\substack{d_1,d_2,\ldots, d_k \le R\\ \mathcal{D}(\mbf{H})}}\left(\prod_{\nu \in \mathcal{P}(k)}\frac{{\mu(a_\nu)}^{|\nu|}}{a_\nu}\right)\left(\prod_{j =1}^k  \log\frac{R}{d_j}\right).\label{4.10}\end{equation}
We now apply the formula, for  $c>0$, 
\begin{equation} \frac{1}{ 2\pi i}\int_{c-i\infty}^{c+i\infty} \frac{x^s}{ s^2}\, ds = 
\left\{ \begin{array}{ll}
       0, 
        &\text{if $0<x\le 1$, 
	}  \\
\log x,& \text{if $x\ge 1$,}   
\end{array}
\right. \label{4.11}
\end{equation}
and have that
\begin{equation}
\mathcal{T}_k(\mathcal{H}_1,\mathcal{H}_2)= \frac {1}{ (2\pi i)^{k}}\mathop{\int }_{(c_k)\ }\! \cdots \!\mathop{\int}_{(c_2)\ } \!\mathop{\int}_{(c_1)\ }F(s_1,s_2,\ldots, s_k)\prod_{j =1}^k\frac{R^{s_j}}{ {s_j}^2} ds_j ,\label{4.12}\end{equation}
where
\begin{equation} F(s_1, s_2, \ldots , s_k) = \sumprime_{\substack{a_\nu,\  \nu\in  \mathcal{P}(k)\\ \mathcal{D}(\mbf{H})  }}\prod_{\nu\in \mathcal{P}(k)}\frac{{\mu(a_\nu)}^{|\nu|}}{{a_\nu}^{1 + \tau_\nu}}, \label{4.13} \end{equation}
and 
\begin{equation} \tau_\nu =\sum_{j\in \nu}s_j .\label{4.14} \end{equation}

We next consider the divisibility conditions $\mathcal{D}(\mbf{H})$. 
The variables $a_\nu$ indexed by the singleton sets $\nu=\{j\}$, $1\le j\le k$ are not constrained by these divisibility conditions, and therefore can contain any prime as a factor.  Further, if $r\ge 1$, then $h_j-h_i=0$ for $j=k-i+1$ and $1\le i\le r$. Thus these constraints drop out of $\mathcal{D}(\mbf{H})$ and the unconstrained variables are both the singleton sets $\nu = \{j\}$, $1\le j\le k$, and also the doubleton sets $\nu=\{i,k-i+1\}$, $1\le i\le r$. (If $r=0$ there are none of these doubleton sets.)  The remaining $a_\nu$ are constrained by at least one of the divisibility relations, and therefore must divide some $h_j-h_i$ so that they can only contain prime factors $\le h$. 
We therefore see that we can write $F(s_1, \ldots , s_k)$ as the Euler product, for $\sigma_j>0$, $1\le j\le k$,          
\begin{equation}\begin{split}  F(s_1, \ldots , s_k)  &= \prod_{p\le h}\Big( 1 -\sum_{j=1}^k \frac{1}{p^{1+ s_j}}  + f_{\mbf{H}}(p;s_1,s_2,\ldots, s_k)\Big)\\& \hskip 1in \times \prod_{p>h}\Big( 1 -\sum_{j=1}^k \frac{1}{p^{1+ s_j}}+\sum_{j=1}^r \frac{1}{p^{1+ s_j+s_{k-j+1}}}\Big),
\label{4.15}\end{split} \end{equation}
where 
\begin{equation} f_{\mbf{H}}(p;s_1,s_2,\ldots ,s_k)=\sum_{\substack{\nu\in {\mathcal{P}}(k),|\nu|\ge 2 \\ p|h_j-h_i\ \textrm{for all } i,j\in \nu } }\frac{(-1)^{|\nu|}}{p^{1+\tau_\nu}}.\label{4.16}\end{equation}
Factoring out the dominant zeta-factors we write
\begin{equation}\begin{split}& F(s_1,s_2,\ldots , s_k)\\& \quad = G_{\mbf{H}}(s_1,s_2,\ldots , s_k) \prod_{j=1}^r\frac{\zeta(1+s_j+s_{k-j+1})}{\zeta(1+s_j)\zeta(1+s_{k-j+1})}\prod_{j=r+1}^\kappa \frac{1}{\zeta(1+s_j)},\label{4.17} \end{split}\end{equation}
and proceed to analyze $G_{\mbf{H}}$.  Let
\begin{equation} \Delta := \prod_{1\le i<j\le \kappa} |h_j-h_i| \le h^{k^2} ,\label{4.18}\end{equation}
so that this product is over all the non-zero differences of  $h_i$ and $h_j$ for $1\le i<j\le k$. 
(Here of course $\Delta$ is not the same function as in the first section.) From the discussion above \eqref{4.15}, $f_{\mbf{H}}=0$ unless $p|\Delta$, and therefore 
\begin{equation} G_{\mbf{H}}(s_1,s_2,\ldots, s_k) = \prod_{p|\Delta}\left(\frac{\displaystyle 1 -\sum_{j=1}^k \frac{1}{p^{1+ s_j}} + f_{\mbf{H}}(p;s_1,s_2,\ldots, s_k)}{\displaystyle \prod_{j=1}^k\Big(1- \frac{1}{p^{1+ s_j}}\Big)\prod_{j=1}^r\Big(1- \frac{1}{p^{1+ s_j+s_{k-j+1}}}\Big)^{-1}}\right)h(s_1,s_2,\ldots s_k), \label{4.19}\end{equation}
where
\begin{equation} h(s_1,s_2,\ldots s_k) = \prod_{p\ndiv\Delta}\left(\frac{\displaystyle 1 -\sum_{j=1}^k \frac{1}{p^{1+ s_j}} +\sum_{j=1}^r \frac{1}{p^{1+ s_j+s_{k-j+1}}}}{\displaystyle \prod_{j=1}^k\Big(1- \frac{1}{p^{1+ s_j}}\Big)\prod_{j=1}^r\Big(1- \frac{1}{p^{1+ s_j+s_{k-j+1}}}\Big)^{-1}}\right). \label{4.20}\end{equation}
Let
\begin{equation} s^*= -\sum_{j=1}^k \min(\sigma_j,0).\label{4.21}\end{equation}
Taking $\sigma_j \ge -\frac{1}{5}$, $1\le j\le k$, we have
\begin{equation}h(s_1,s_2,\ldots s_k) \ll_k \prod_p\Big(1 +O_k(\frac{1}{p^{6/5}})\Big) \ \ll_k 1, \label{4.22} \end{equation}
and thus in this region we have
\begin{equation} \begin{split}G_{\mbf{H}}(s_1,s_2,\ldots , s_k)&\ll_k\prod_{p|\Delta}\Big(1+O_k(\frac{1}{p^{1-s^*}})\Big) \\& \ll_k \exp\left(a(k) \sum_{p\le {k^2}\log(2h)}\frac{1}{p^{1-s^*}}\right) \\& \ll_k \exp\left(a(k)({k^2}\log(2h))^{s^*} \sum_{p\le {k^2}\log(2h)}\frac{1}{p}\right) \\& \ll_k 
       \exp\left(b(k)(\log(2h))^{s^*}\log\log\log 2h\right),
\label{4.23} \end{split}\end{equation}
where the sum which was originally over $p|\Delta$ has been majorized by using \eqref{4.18} to find the smallest set of primes that could divide $\Delta$.
By this bound and \eqref{4.17} we see that if $r\ge 1$ then $F$  has  simple poles at $s_i+s_{k-i+1}=0$, $1\le i\le r$. By \eqref{3.4}, for $|s_i|\ge C/2$, 
\begin{equation} F(s_1,s_2,\ldots , s_k)   \ll_k \exp( b(k)(\log(2h))^{s^*}\log \log\log 2h)\prod_{j=1}^k\log^2(|t_j|+2)\prod_{i=1}^r\frac{1}{|s_i+s_{k-i+1}|}
.\label{4.24}\end{equation}

We are now ready to begin the evaluation of $\mathcal{T}_k(\mathcal{H}_1,\mathcal{H}_2)$.  By \eqref{4.24} we see the integrand  in \eqref{4.12} goes to zero as any one of the variables $|t_j|\to \infty$, and thus we can move any contours we wish to the left to $\mathcal{L}$. We first move successively the contours $(c_{j})$, $r+1\le j \le  \kappa$ to  $\mathcal{L}$; by \eqref{4.17} these correspond to the cases where the integrand has only a simple pole at $s_j=0$. If $r=\kappa$ there are none of these terms and we skip ahead to \eqref{4.28}. Thus, moving $c_{r+1}$ to $\mathcal{L}$  and passing a simple pole at $s_{r+1}=0$ we obtain
\begin{equation}\begin{split}
&\mathcal{T}_k(\mathcal{H}_1,\mathcal{H}_2)\\& = \frac {1}{ (2\pi i)^{k-1}}\bigg(\prod_{\substack{ j=1\\j\neq r+1}}^k\mathop{\int }_{(c_{j})\ }\!\bigg) G_{\mbf{H}}(s_1,s_2,\ldots, s_{k})\Bigg|_{s_{r+1}=0}\bigg(\prod_{ j=r+2}^\kappa \frac{R^{s_j}}{\zeta(1+s_j){s_j}^2} ds_j\bigg) \\& \hskip 1in \times \bigg(\prod_{j=1}^{r}\frac{\zeta(1+s_j+s_{k-j+1})R^{s_j+s_{k-j+1}}}{\zeta(1+s_j)\zeta(1+s_{k-j+1}){s_j}^2(s_{k-j+1})^2}ds_j\, ds_{k-j+1}\bigg)
\\&\hskip 1.2in  + \frac {1}{ (2\pi i)^{k}}\bigg(\prod_{\substack{ j=1\\j\neq r+1}}^k\mathop{\int }_{(c_{j})\ }\!\bigg) \mathop{\int}_{\mathcal{L}\ }F(s_1,s_2,\ldots, s_{k})\prod_{j =1}^k\frac{R^{s_j}}{ {s_j}^2} ds_j.\label{4.25}\end{split}\end{equation}
We bound the second multiple integral on the right by moving all the contours $(c_j)$, $j\neq r+1$, to $(\frac{1}{\log R})$ which leaves the value of the integral unchanged. If $s_j$ and $s_{k-j+1}$ are on $(\frac{1}{\log R})$, 
\begin{equation} \begin{split}&\frac{\zeta(1+s_j+s_{k-j+1})}{\zeta(1+s_j)\zeta(1+s_{k-j+1})}\\&\hskip .5in \ll (\log R) \log(2+|s_j+s_{k-j+1}|)\log(2+|s_j|)
\log(2+|s_{k-j+1}|). \end{split}\label{4.26}\end{equation}
 In the multiple integral   $s^*= -\sigma_{r+1} \le \frac{C}{\log 2}$ for $\sigma_{r+1}$ on $\mathcal{L}$,  and  we  take a fixed $C< \frac{\log 2}{2}$ . Then by \eqref{4.17}, \eqref{4.24}, \eqref{4.26}  and Lemma 1 we have the second term in \eqref{4.25} is
\begin{equation} \begin{split} &\ll_k \exp\left( b(k)(\log(2h))^{\frac{1}{2}}\right)\bigg(\mathop{\int}_{(\frac{1}{\log R})\ }\log^2(|s|+2)\left|\frac{R^s ds}{s^2}\right|\bigg)^{\kappa -r-1}(\log R)^r\\& \hskip 1in \times \bigg(\mathop{\int}_{(\frac{1}{\log R})\ }\log^2(|s|+2)\left|\frac{R^s ds}{s^2}\right|\bigg)^{2r}\times \mathop{\int}_{\mathcal{L}\ }\log^2(|s|+2)\left|  \frac {R^s ds}{s^2}  \right| \\&\ll_k  \exp\left( b(k)(\log(2h))^{\frac{1}{ 2}}\right)(\log R)^{\kappa +2r-1}e^{-c'\sqrt{\log R}}
\\&
\ll_k e^{-{c'}_k\sqrt{\log R}},\label{4.27} \end{split}\end{equation}
where we used $\log 2h \ll \log R$ for the last line.

We continue this process, moving next $(c_{r+2})$ to $\mathcal{L}$ in the first term, and estimating the secondary term as above, and so on successively through the contours up to $(c_{\kappa})$. Hence we conclude 
\begin{equation}  \mathcal{T}_k(\mathcal{H}_1,\mathcal{H}_2) 
=\mathcal{U}_r +O_k(e^{-{c'}_k\sqrt{\log R}}), \label{4.28} \end{equation} 
where
\begin{equation}\begin{split} \mathcal{U}_r &= \frac {1}{ (2\pi i)^{2r}}\bigg(\prod_{j=1}^r\mathop{\int }_{(c_{k-j+1})\ } \!\mathop{\int}_{(c_{j})\ }\bigg) G_1(s_1,s_2,\ldots , s_r,s_{\kappa+1},s_{\kappa+2},\ldots s_k) \\& \hskip .8in \times\prod_{j=1}^{r}\frac{\zeta(1+s_{j}+s_{k-j+1})R^{s_{j}+s_{k-j+1}}}{\zeta(1+s_{j})\zeta(1+s_{k-j+1}){s_{j}}^2(s_{k-j+1})^2}ds_{j}\, ds_{k-j+1},\end{split}\label{4.29}\end{equation}
and
\begin{equation} G_1(s_1,s_2,\ldots , s_r,s_{\kappa+1},s_{\kappa+2},\ldots , s_k)= G_{\mbf{H}}(s_1,s_2,\ldots , s_{r},0,0,\ldots , 0,s_{\kappa+1},s_{\kappa+2}, \ldots , s_k). \label{4.30} \end{equation}

We will now prove that 
\begin{equation} \mathcal{U}_r = G_{\mbf{H}}(0,0,\ldots, 0)(\log R)^r + \sum_{j=1}^r\mathcal{A}_j(\mathcal{H})(\log R)^{r-j} + O_k(e^{-{c'}_k\sqrt{\log R}}), \label{4.31}\end{equation}
where the $\mathcal{A}_j(\mathcal{H})$ are explicitly computable arithmetic functions which for $1\le h \le R^A$ with any $A>0$ satisfy the bound
\begin{equation}\mathcal{A}_j(\mathcal{H}) \ll_k (\log\log 2h)^{b(k)}. \label{4.32}\end{equation}
We will also prove at the end of this section that
\begin{equation}  G_{\mbf{H}}(0,0,\ldots, 0) = \gs(\mathcal{H}). \label{4.33} \end{equation}
From these results Proposition 1 follows.  

The multiple integral in $\mathcal{U}_r$ would decouple into a product of double integrals evaluated in Lemma 2 if $G_1$ were a constant, but since this is not the case we need to apply Lemma 2 inductively. To do this we need estimates for the partial derivatives of $G_{\mbf{H}}$. Let $\mbf{a}=(a_1,a_2,\ldots , a_k)$, and define
\begin{equation} D_{\mbf{a}}G_{\mbf{H}} = \frac{\partial ^{a_1}}{\partial {s_1}^{a_1}}\frac{\partial ^{a_2}}{\partial {s_2}^{a_2}}\cdots \frac{\partial ^{a_k}}{\partial {s_k}^{a_k}}\, G_{\mbf{H}}(s_1,s_2,\ldots ,s_k) .\label{4.34} \end{equation}
We have, for $\sigma_j > -{c'}_k $, $1\le j \le k$, 
\begin{equation} D_{\mbf{a}}G_{\mbf{H}}\ll_k
       (\log \log 2h)^{b'(k)}\exp\left( b(k)(\log(2h))^{s^*}\log\log\log 2h\right). 
\label{4.35} \end{equation}
To obtain these estimates, we logarithmically differentiate $G_{\mbf{H}}$ to see 
\[ \frac{\partial G_{\mbf{H}}}{\partial s_1} \ll_k |G_{\mbf{H}}|\left(\sum_{p\le k^2\log(2h)} \frac{\log p}{p^{1-s^*}}\right) ,\]
The sum above is bounded by
\[\ll (k^2\log(2h))^{s^*} \sum_{p\le k^2\log(2h)} \frac{\log p}{p} \ll (k^2\log(2h))^{s^*}  \log (k^2\log(2h)), \] 
and \eqref{4.35} follows in this case by \eqref{4.23}.
 By the product rule, further partial derivatives will satisfy the above bound with the sum having $\log p$ replaced   by $(\log p)^{c(k)}$, which only changes the value of $b'(k)$ in \eqref{4.35}.

We first consider the case $r=1$ in \eqref{4.29}. By Lemma 2 applied with $f_R=G_1$ we see by \eqref{4.23} that the conditions of the lemma are satisfied and therefore
\[ \mathcal{U}_1 =   G_{\mbf{H}}(0,0,\ldots,0)\log R + \mathcal{A}_1 + O_k(e^{-{c'}_k\sqrt{\log R}}), \]
where
\[ \mathcal{A}_1 =  \frac{\partial G_1}{\partial s_k}(0,0,\ldots , 0) + \frac {1}{ 2\pi i}\mathop{\int}_{\mathcal{L}\ } G_1(s_1,-s_1) \frac{ds_1}{\zeta(1+s_1)\zeta(1-s_{1}){s_1}^4}. \]
It remains to prove that $\mathcal{A}_1$ satisfies \eqref{4.32}. By \eqref{4.35} the first term in $\mathcal{A}_1$ satisfies this bound. In the integral term we move the contour back to the imaginary axis with a semicircle of radius $\delta=\frac{1}{\log(k^2\log(2h))}$ to the left of the double pole at $s_1=0$.  Using \eqref{3.4} and \eqref{4.23} the part of the integral over the contour on the imaginary axis is bounded by
\[ \begin{split}&\ll_k  (\log \log 2h )^{b(k)} \int_{\delta}^\infty \frac{(\log (t+2))^2}{t^4}\, dt \\& \ll_k\frac{ (\log \log 2h )^{b(k)}}{\delta^3}\\&
\ll_k  (\log \log 2h )^{b'(k)} , \end{split}\]
and the integral over the contour on the semicircle is bounded by
\[ \begin{split}  &\ll_k \pi  \delta \times \frac{ (\log \log 2h )^{3b(k)}}{\delta^2}
\\& \ll_k  (\log \log 2h )^{b'(k)} ,\end{split}\]
which completes the proof for $r=1$.

 For the general case of \eqref{4.29}, we move all the contours to $(\frac{1}{\log R})$ and apply Lemma 2  for the double integral over $s_1$ and $s_k$ to obtain
\begin{equation}\begin{split} \mathcal{U}_r &= \frac {1}{ (2\pi i)^{2r-2}}\mathop{\int }_{(\frac{1}{\log R})\ }\! \cdots \!\mathop{\int}_{(\frac{1}{\log R})\ }\bigg( G_1\bigg|_{\substack{s_1=0\\s_k=0}}\log R +G_2+ O_r(e^{-c_r'\sqrt{\log R}})\bigg)\\&\hskip 1.1in \times \prod_{j=2}^{r}\frac{\zeta(1+s_{j}+s_{k-j+1})R^{s_{j}+s_{k-j+1}}}{\zeta(1+s_j)\zeta(1+s_{k-j+1}) (s_{k-j+1})^2{s_{j}}^2} ds_{j} ds_{k-j+1}\\& =\mathcal{U}_{r-1}\log R + {\mathcal{U}'}_{r-1}+ O_r(e^{-c_r'\sqrt{\log R}}),\label{4.36}\end{split}\end{equation}
where the error term was estimated using \eqref{4.26} as in \eqref{4.27}. Here 
\[ \begin{split} G_2(s_2,s_3,\ldots,& s_r,s_{\kappa+1},s_{\kappa +2}, \ldots, s_{k-1}) = 
\frac{\partial G_1}{\partial s_k}(0,s_2,s_3,\ldots, s_r,s_{\kappa+1},s_{\kappa +2}, \ldots, s_{k-1},0)
\\& + \frac {1}{ 2\pi i}\mathop{\int}_{\mathcal{L}\ } G_1(s_1,s_2,\ldots, s_r,s_{\kappa+1}, \ldots , s_{k-1}, -s_1) \frac{ds_1}{\zeta(1+s_1)\zeta(1-s_{1}){s_1}^4}, \end{split}\]
and therefore $\mathcal{U'}_{r-1}$ is of the same form as $\mathcal{U}_{r-1}$ with $G_1$ replaced by a partial derivative of $G_1$ or an absolutely convergent integral of $G_1$ with respect to the variable $s_1$ when $s_k=-s_1$. As we saw in the case $r=1$, both of these terms satisfy \eqref{4.32} and \eqref{4.35}. We now apply Lemma 2 for the double integral over $s_2$ and $s_{k-1}$, and continue this process until all the variables are exhausted. We thus arrive at \eqref{4.31} and the bound \eqref{4.32} follows by \eqref{4.35} and the argument used in the case $r=1$.

It remains to prove \eqref{4.33}. By \eqref{4.19} and \eqref{4.20} we have  
\[ G_{\mbf{H}}(0,0,\ldots, 0) = \prod_{p|\Delta} \left(1 -\frac{k}{p} +f_{\mbf{H}}(p;0,0,\ldots,0)\right)\left(1-\frac{1}{p}\right)^{-\kappa} \prod_{p\ndiv \Delta} \left(1 -\frac{\kappa}{p} \right)\left(1-\frac{1}{p}\right)^{-\kappa},\]
where by \eqref{4.16}
\[f_{\mbf{H}}(p;0,0,\ldots, 0)= \frac{1}{p}\sum_{\substack{\nu\in {\mathcal{P}}(k),|\nu|\ge 2 \\ p|h_j-h_i\ \textrm{for all } i,j\in \nu }}  (-1)^{|\nu|}. \]
Therefore by \eqref{2.2} we need to prove
\begin{equation}\sum_{\substack{\nu\in {\mathcal{P}}(k),|\nu|\ge 2 \\ p|h_j-h_i\ \textrm{for all } i,j\in \nu }}  (-1)^{|\nu|}= k-\nu_p(\mathcal{H}). \label{4.37}\end{equation}
If $\nu_p(\mathcal{H})=q$, then $h_1, h_2,\ldots , h_k$ must fall into $q$ distinct residue classes, say $r_1,r_2,\ldots, r_q$ $(\mathrm{mod}\ p)$. Let
\[ \mathcal{M}_p(\ell) = \{m: h_m\equiv r_\ell (\mathrm{mod}\ p)\}, \quad 1\le \ell \le q .\] 
Thus given $p$ the sets $\mathcal{M}_p(\ell)$, $1\le \ell\le q$, give a disjoint partition of the set $\{1,2,\ldots, k\}$, and therefore
\begin{equation} \sum_{\ell=1}^q |\mathcal{M}_p(\ell)| =k.\label{4.38}\end{equation}
The conditions $p|h_j-h_i$ hold if and only if $h_i$ and $h_j$ are in the same residue class modulo $p$ and thus if and only if $i$ and $j$ are in $\mathcal{M}_p(\ell)$ for some $\ell$. 
Hence the $\nu \in \mathcal{P}(k)$ which will satisfy $p|h_j-h_i$ for all $i,j\in \nu$ are precisely the subsets of $\mathcal{M}_p(\ell)$  with at least two elements 
\[ \tilde{\mathcal{M}}_p = \bigcup_{\ell=1}^q \{\nu: \nu\subset \mathcal{M}_p(\ell),|\nu|\ge 2\}.\]
We conclude, using \eqref{4.38}, that
\[\begin{split} \sum_{\substack{\nu\in {\mathcal{P}}(k),|\nu|\ge 2 \\ p|h_j-h_i\ \textrm{for all } i,j\in \nu }}  (-1)^{|\nu|}&=\sum_{\nu \in \tilde{\mathcal{M}}_p }(-1)^{|\nu|}\\  & = \sum_{\ell=1}^q \sum_{\nu \subset \mathcal{M}_p(\ell), |\nu|\ge 2 }(-1)^{|\nu|}\\  &=
\sum_{\ell=1}^q \sum_{j=2}^{|\mathcal{M}_p(\ell)|}(-1)^j\binom{|\mathcal{M}_p(\ell)|}{j} \\
&= \sum_{\ell=1}^q\Big( -1 + |\mathcal{M}_p(\ell)|\Big)   \\
&=  \sum_{\ell=1}^q |\mathcal{M}_p(\ell)| -q\\ 
& = k - \nu_p(\mathcal{H}).\end{split}\]

\section{Proof of Proposition 2}
We first reduce the proof to the special case when $h_0 \not \in \mathcal{H}=\mathcal{H}_1\cup \mathcal{H}_2$.
Let
\begin{equation} \tilde{\mathcal{S}}(\mathcal{H}_1,\mathcal{H}_2,h_0) = \sum_{n\le N} \Lambda_R(n;\mathcal{H}_1)\Lambda_R(n;\mathcal{H}_2)\Lambda(n+h_0) . \label{5.1} \end{equation}
Since trivially $|\Lambda_R(n)|\le d(n)\log R$, we see for $i=1,2$
\begin{equation} \Lambda_R(n,\mathcal{H}_i) \le (d(n)\log R)^{k_i},  \label{5.2} \end{equation}
and since $d(n)\ll n^\epsilon$ and in Proposition 2 $R\le N^{\frac{1}{2}}$, we have
\begin{equation} \begin{split} \tilde{\mathcal{S}}(\mathcal{H}_1,\mathcal{H}_2,h_0) &= \sum_{R<n\le N} \Lambda_R(n;\mathcal{H}_1)\Lambda_R(n;\mathcal{H}_2)\Lambda(n+h_0) +O(R^{1+\epsilon}) \\
&=\sum_{\substack{R<n\le N\\ n+h_0\ \mathrm{prime}}} \Lambda_R(n;\mathcal{H}_1)\Lambda_R(n;\mathcal{H}_2)\Lambda(n+h_0) +O(N^{\frac{1}{2}+\epsilon}) , \label{5.3} \end{split}\end{equation}
where we have removed the prime powers in the last line. If  $n+h_0$ is a prime $>R$ then its only divisor $\le R$ is $d=1$, and therefore
\[ \Lambda_R(n+h_0)\Lambda(n+h_0) =  \Lambda(n+h_0)\log R. \]
Thus if  $h_0 \in  \mathcal{H}_1 \cap \mathcal{H}_2 $,
\[\begin{split} &\tilde{\mathcal{S}}(\mathcal{H}_1,\mathcal{H}_2,h_0) 
\\& =(\log R)^2 \sum_{\substack{R<n\le N\\ n+h_0\ \mathrm{prime}}} \Lambda_R(n;\mathcal{H}_1-\{h_0\})\Lambda_R(n;\mathcal{H}_2-\{h_0\})\log(n+h_0) +O(N^{\frac{1}{2}+\epsilon});\end{split}\]
if $h_0\in \mathcal{H}_1 - \mathcal{H}_2$,
\[ \begin{split}\tilde{\mathcal{S}}(\mathcal{H}_1,\mathcal{H}_2,h_0) 
= \log R \sum_{\substack{R<n\le N\\ n+h_0\ \mathrm{prime}}} \Lambda_R(n;\mathcal{H}_1-\{h_0\})\Lambda_R(n;\mathcal{H}_2)&\log(n+h_0) \\& +O(N^{\frac{1}{2}+\epsilon});\end{split} \]
and if  $h_0 \not \in \mathcal{H}_1 \cup \mathcal{H}_2 $, 
\[ \tilde{\mathcal{S}}(\mathcal{H}_1,\mathcal{H}_2,h_0) 
=  \sum_{\substack{R<n\le N\\ n+h_0\ \mathrm{prime}}} \Lambda_R(n;\mathcal{H}_1)\Lambda_R(n;\mathcal{H}_2)\log(n+h_0) +O(N^{\frac{1}{2}+\epsilon}).  \]
In these sums we may once again include the terms $\le R$ and the prime powers if we wish with the same error term, and therefore in each situation we have reduced the proof to the case when $h_0$ is distinct from the other $h_i$'s. Henceforth we can therefore take 
\begin{equation} h_0 \not \in \mathcal{H}.\label{5.4}\end{equation}

Proceeding to the proof, we have 
\begin{equation} \tilde{\mathcal{S}}(\mathcal{H}_1,\mathcal{H}_2,h_0)  =\sum_{d_1,d_2,\ldots , d_k\le R}\Big(\prod_{i=1}^k \mu(d_i)\log \frac{R}{d_i}\Big) \sum_{\substack{n\le N\\ d_i|n+h_i,\ 1\le i\le k }}\Lambda(n+h_0) . \label{5.5} \end{equation}
By the Chinese Remainder Theorem the sum will run through an arithmetic progression modulo $D_k$  provided $(d_i,d_j)|h_j-h_i$, $1\le i<j\le k$ and will be empty otherwise. As in \eqref{4.9} we denote these conditions by $\mathcal{D}(\mathcal{H})$. Using Iverson notation, we let 
\begin{equation} \psi(x;q,a) := \sum_{\substack{n\le x\\ n\equiv a (\textrm{mod}\ q) }}\Lambda(n) = [(a,q)=1]\frac{x}{\phi(q)} + E(x;q,a), \label{5.6} \end{equation} 
and have
\begin{equation} \begin{split} \sum_{\substack{n\le N\\ d_i|n+h_i,\ 1\le i\le k }}\Lambda(n+h_0) &= [\mathcal{D}(\mathcal{H})]\Big(\psi(N+h_0;D_k,a) - \psi(h_0,D_k,a)\Big)\\&= [\mathcal{D}(\mathcal{H})]\psi(N;D_k,a) + O(h\log N), \label{5.7}\end{split} \end{equation}
where $a$ is an integer satisfying the congruence relations $a\equiv h_0-h_j (\textrm{mod}\ d_j)$, $1\le j\le k$. The term $\psi(N;D_k,a)$ has a non-zero main term if $(a,D_k)=1$, which is equivalent to
\begin{equation}  (d_j, h_j-h_0)=1, \ 1\le j\le k,\label{5.8} \end{equation}
and, if $(a,D_k)>1$ then $\psi(N;D_k,a)\ll (\log N)^2$; thus
\begin{equation} \begin{split} \tilde{\mathcal{S}}(\mathcal{H}_1,&\mathcal{H}_2,h_0) =N \sum_{\substack{d_1,d_2,\ldots d_k\le R\\ \mathcal{D}(\mathcal{H})\\(d_j, h_j-h_0)=1, \ 1\le j\le k }}\frac{1}{\phi(D_k)}\Big(\prod_{i=1}^k \mu(d_i)\log \frac{R}{d_i}\Big)\\&\qquad + O\big(\sum_{d_1,d_2,\ldots d_k\le R}\Big(\prod_{i=1}^k \mu^2(d_i)\log \frac{R}{d_i}\Big)\max_{\substack{a(\mathrm{mod}\ D_k)\\ (a,D_k)=1}}|E(N;D_k,a)|\big) \\ & \qquad  +O\big( R^k h(\log N)^2\big) \\& = N\tilde{\mathcal{T}}_k(\mathcal{H}_1,\mathcal{H}_2,h_0) +O( \mathcal{E}_k) + O\big( R^kh(\log N)^2\big)  . \label{5.9} \end{split}\end{equation}
We handle the error term $\mathcal{E}_k$ with the Bombieri-Vinogradov Theorem. First, 
\[\begin{split} \mathcal{E}_k &\ll (\log R)^k \sum_{d_1,d_2,\ldots , d_k\le R}\mu^2(D_k)\max_{\substack{a(\mathrm{mod}\ D_k)\\ (a,D_k)=1}}|E(N;D_k,a)| \\& \ll  (\log R)^k \sum_{q\le R^k}\mu^2(q)\max_{\substack{a(\mathrm{mod}\ q)\\ (a,q)=1}}|E(N;q,a)|\sum_{\substack{q=D_k\\ d_1,d_2,\ldots ,d_k\le R}}1.
 \end{split}\]
Given $q$, the number of ways to write $q=D_k$ (i.e. write $q$ as the least common multiple of $k$ squarefree numbers) is bounded by $d(q)^k$, since each of the $k$ numbers in the least common multiple must be a divisor of $q$.  Applying Cauchy's inequality, we have
\[\begin{split} \mathcal{E}_k &\ll (\log R)^k \sum_{q\le R^k}\mu^2(q)d(q)^k \max_{\substack{a(\mathrm{mod}\ q)\\ (a,q)=1}}|E(N;q,a)|\\& \ll (\log R)^k \sqrt{ \sum_{q\le R^k}\frac{{d(q)}^{2k}}{q}}\sqrt{\sum_{q\le R^k}q \max_{\substack{a(\mathrm{mod}\ q)\\ (a,q)=1}}|E(N;q,a)|^2}.\end{split} \]
We now use the estimate
\begin{equation} \sum_{n\le x} d(n)^k \ll_k x (\log x)^{2^k-1} \label{5.10}\end{equation}
and the trivial estimate $|E(N;q,a)|\ll \frac{N\log N}{q}$ for $q\le N$  to conclude
\[ \mathcal{E}_k \ll_k (\log R)^{4^k +k }\sqrt{N\log N} \sqrt{\sum_{q\le R^k} \max_{\substack{a(\mathrm{mod}\ q)\\ (a,q)=1}}|E(N;q,a)|}.\]
By the Bombieri-Vinogradov Theorem the sum on is $\ll \frac{N}{(\log N)^A}$ for any $A>0$ provided
\begin{equation}R^k \le N^{\frac{1}{2}}(\log N)^{-B},\label{5.11} \end{equation}
where $B=B(A)$.  We conclude under this condition that
\begin{equation} \mathcal{E}_k\ll_k N(\log N)^{4^k+k+\frac{1}{2}-\frac{A}{2}}=o_k(N) \label{5.12}\end{equation}
if $A> 2(4^k+k+\frac{1}{2})$. 
To complete the proof of the proposition we will now prove that, for $R^k \le N$ and $h\le R^{\frac{1}{2k}}$,
\begin{equation} \tilde{\mathcal{T}}_k(\mathcal{H}_1,\mathcal{H}_2,h_0) = \mathcal{T}_{k+1}(\mathcal{H}_1\cup \{h_0\},\mathcal{H}_2) +O_k(e^{-c_k\sqrt{\log R}}), \label{5.13} \end{equation}
which by \eqref{4.28}, \eqref{4.32}, and \eqref{4.33} completes the proof. 
To prove \eqref{5.13}, we have
\begin{equation}\begin{split}&\mathcal{T}_{k+1}(\mathcal{H}_1\cup \{h_0\},\mathcal{H}_2) = \sum_{\substack{d_0,d_1,\cdots, d_k\le R\\(d_i,d_j)|h_j-h_i,0\le i<j\le k}} \frac{1}{[d_0,D_k]}\prod_{j=0}^k\mu(d_j)\log \frac{R}{d_j}\\&= \sum_{\substack{d_1,d_2,\cdots, d_k\le R\\(d_i,d_j)|h_j-h_i,\ 1\le i<j\le k}}\Big( \prod_{j=1}^k\mu(d_j)\log \frac{R}{d_j}\Big)\sum_{\substack{d_0\le R\\ (d_0,d_j)|h_j-h_0,\ 1\le j\le k}}\frac{\mu(d_0)}{[d_0,D_k]}\log\frac{R}{d_0} \\ &
=\sum_{\substack{d_1,d_2,\cdots, d_k\le R\\\mathcal{D}(\mathcal{H})}}\Big( \prod_{j=1}^k\mu(d_j)\log \frac{R}{d_j}\Big)T_{1}(\mathcal{H}_1\cup \{h_0\},\mathcal{H}_2).\end{split}\label{5.14}\end{equation}
On letting $g=(d_0,D_k)$, $d_0=gd'$, we see $[d_0,D_k]=D_k d'$ and $(d',D_k)=1$. Thus 
\[ \begin{split} T_{1}(\mathcal{H}_1\cup \{h_0\},\mathcal{H}_2)&= \sum_{\substack{ gd'\le R\\ g|D_k\\ (g,d_j)|h_j-h_0, \ 1\le j\le k\\ (d',D_k)=1}} \frac{\mu(g d')}{d'D_k}\log\frac{R}{g d'} \\
&=\frac{1}{D_k} \sum_{\substack{g\le R\\ g|D_k\\ (g,d_j)|h_j-h_0, \ 1\le j\le k}} \mu(g) \sum_{\substack{d'\le R/g \\ (d',D_k)=1}}\frac{\mu( d')}{d'}\log\frac{R/g}{ d'} .\end{split}\]
For $\log m \ll \log R$ we have (by the prime number theorem or see Lemma 2.1 of \cite{GYI})
\begin{equation} \sum_{\substack{d\le R\\ (d,m)=1}} \frac{\mu(d)}{d}\log \frac{R}{d} = \frac{m}{\phi(m)} + O(e^{-c_1\sqrt{\log R}}). \label{5.15} \end{equation}
Applying this and dropping the redundant condition $g\le R$  since $g\le \prod_{j=1}^k (g,d_j) \le h^k \le R $ when $h\le R^{1/k}$, we see 
\[ T_{1}(\mathcal{H}_1\cup \{h_0\},\mathcal{H}_2) = \frac{1}{\phi(D_k)} \sum_{\substack{g|D_k\\ (g,d_j)|h_j-h_0,\ 1\le j\le k}} \mu(g) + O(\frac{d(D_k)}{D_k}e^{-c_1\sqrt{\log \left(R/h^k\right)}}).\]
We now claim that
\[\sum_{\substack{g|D_k\\ (g,d_j)|h_j-h_0,\ 1\le j\le k}} \mu(g)= [(d_j,h_j-h_0)=1 ,\ 1 \le j\le k].\]
One way to see this is through the decomposition of the $d_i$'s into relatively factors \eqref{4.8} from which we see we can write $g=\prod_{\nu \in \mathcal{P}(k)}b_\nu$, $b_\nu |a_\nu$ with the $b_\nu$ pairwise relatively prime with each other. Then the sum becomes
\[ \prod_{j=1}^k \prod_{\nu \in \mathcal{P}_{\{j\}}(k)}\sum_{b_\nu|(d_j,h_j-h_0)}\mu(b_\nu)=[(d_j,h_j-h_0)=1 ,\ 1 \le j\le k] .\]
We conclude
\[\begin{split} T_{1}(\mathcal{H}_1\cup \{h_0\},\mathcal{H}_2)=  \frac{1}{\phi(D_k)}&[(d_j,h_j-h_0)=1 ,\ 1 \le j\le k]\\& + O_k(\frac{d(D_k)}{D_K}e^{-c_1\sqrt{\log \left(R/h^k\right)}}),\end{split}\]
and on substituting this result in \eqref{5.14} we have 
\[\begin{split} \mathcal{T}_{k+1}&(\mathcal{H}_1\cup \{h_0\},\mathcal{H}_2) = \sum_{\substack{d_1,d_2,\cdots d_k\le R\\\mathcal{D}(k)\\(d_j,h_j-h_0)=1 ,\ 1 \le j\le k}}\frac{1}{\phi(D_k)}\Big( \prod_{j=1}^k\mu(d_j)\log \frac{R}{d_j}\Big)\\& \hskip 1in + O_k((\log R)^k e^{-c_1\sqrt{\log R/h^k}}\sum_{d_1,d_2,\ldots d_k\le R} \frac{d(D_k)}{D_k}).\end{split} \]
The first term is $\tilde{\mathcal{T}}_k(\mathcal{H}_1,\mathcal{H}_2,h_0)$ and by \eqref{5.10} 
\[\begin{split} \sum_{d_1,d_2,\ldots d_k\le R} \frac{d(D_k)}{D_k} &\ll \sum_{q\le R^k}\frac{d(q)}{q} \sum_{q=D_k}1 \\ & \ll \sum_{q\le R^k}\frac{d(q)^{k+1}}{q} \\&
\ll_k (\log R)^{2^{k+1}} .\end{split} \]
Thus the error term is 
 \[ \begin{split} & \ll_k (\log R)^{2^{k+1}+k} e^{-c_1\sqrt{\log \left( R/h^k\right)}} \\ & \ll_k e^{-c_k\sqrt{\log R}} ,\end{split}\]
 which proves \eqref{5.13} if $h\le R^{\frac{1}{2k}}$.

\section{Optimization of a Quadratic Form Related to the Poisson Distribution.}
The content of this section and the proof given here was provided to us by  E. Bombieri and P. Deift. The final tool we need for our proof of Theorem 1 is an optimization procedure related to the Poisson distribution.  Let $ X$ be a Poisson random variable with expected value $\lambda$, defined by the discrete p.d.f. 
\begin{equation} p(j) = \text{Prob.}(X=j) = e^{-\lambda}\frac{\lambda^j}{j!} , \quad j=0,1,2,\ldots\ . \label{6.1}\end{equation}
We define an inner product with respect to this density function by
\begin{equation} \begin{split} \langle f(x),g(x)\rangle &= \sum_{j=0}^\infty f(j)g(j)p(j) \\&
= e^{-\lambda}\sum_{j=0}^\infty f(j)g(j)\frac{\lambda^j}{j!} . \label{6.2} \end{split}\end{equation} 
The $k$-th moment of the Poisson distribution is defined by
\begin{equation} \begin{split}\mu_k(\lambda)= E(x^k) = \langle x^k,1\rangle &= e^{-\lambda}\sum_{j=0}^\infty \frac{j^k\lambda^j}{j!} \\&
=e^{-\lambda}\left(\lambda \frac{d}{d\lambda}\right)^k \sum_{j=0}^\infty \frac{\lambda^j}{j!} \\&
=e^{-\lambda}\left(\lambda \frac{d}{d\lambda}\right)^k e^\lambda . \label{6.3} \end{split} \end{equation}
More explicitly, we have 
\begin{equation}\mu_k(\lambda) = \sum_{\nu =1}^k 
 \left\{		\begin{array}{c} k\\ \nu 	
\end{array}
	\right\}
\lambda^\nu ,\label{6.4}\end{equation}
where  $ \left\{\begin{array}{c} k\\ \nu  	
\end{array} \right\} $  denote  the Stirling numbers of the second type, defined to be the number of ways to partition a $k$-set\footnote{A $k$-set denotes a set with $k$ elements.} into $\nu$ non-empty subsets (not counting the order of the subsets).  It is easy to see that
\begin{equation} \left\{\begin{array}{c} k\\ \nu  	
\end{array} \right\}= \nu \left\{\begin{array}{c} k-1\\ \nu  	
\end{array} \right\}+\left\{\begin{array}{c} k-1\\ \nu -1 	
\end{array} \right\} \label{6.5} \end{equation}
since the last element in our $k$-set either is put into its own singleton set or else it is put into one of the $\nu$ subsets which contain some of the earlier elements.  To prove \eqref{6.4} we use the identity 
\begin{equation} \sum_{\nu=0}^j \nu!  \left\{\begin{array}{c} k\\ \nu  	
\end{array} \right\}  \left(\begin{array}{c} j\\ \nu  	
\end{array} \right) =j^k . \label{6.6} \end{equation}
This identity arises from counting the number of partitions of a $k$-set into less than or equal to $j$ sets, where the order of these sets is counted. On one hand there are $j$ choices for where to place each of the $k$ elements, so this number is $j^k$, while on the other hand, if $\nu$ of these $j$ sets are non-empty, then there are $\nu! \left\{\begin{array}{c} k\\ \nu  	
\end{array} \right\} $ such partitions and $\left(\begin{array}{c} j\\ \nu  	
\end{array} \right)$ ways to choose the $\nu$ non-empty sets. Rewriting \eqref{6.6} in the form
\[  \frac{j^k}{j!}=\sum_{\nu=0}^j  \left\{\begin{array}{c} k\\ \nu  	
\end{array} \right\}  \frac{1}{(j-\nu)!}, \]
 multiplying by $\lambda^j e^{-\lambda}$ and summing over $j$, we obtain by \eqref{6.3}
\[ \begin{split}\mu_k(\lambda) &= e^{-\lambda}\sum_{j=0}^\infty \frac{j^k\lambda^j}{j!} \\&=e^{-\lambda}\sum_{j=0}^\infty \lambda^j \sum_{\nu=0}^j  \left\{\begin{array}{c} k\\ \nu \end{array} \right\}  \frac{1}{(j-\nu)!}\\& = \sum_{\nu=0}^k\left\{\begin{array}{c} k\\ \nu  \end{array} \right\} \lambda^\nu ,\end{split}\]
by interchanging the $j$ and $\nu$ summations, which proves \eqref{6.4}.

Our method for finding small gaps between primes leads us to define a second bilinear form
given by 
\begin{equation} \begin{split} \langle f(x),g(x)\rangle_\rho &= \langle x-\rho,f(x)g(x)\rangle\\& = \sum_{j=0}^\infty (j-\rho)f(j)g(j)p(j)  \label{6.7} \end{split}\end{equation} 
where $\rho$ is a real number. (This is not an inner product because it is not necessarily non-negative.)
Letting $\mathbf{a} = (a_0,a_1,a_2,\ldots , a_k)$, consider
\begin{equation} P_{\mathbf{a}}(x) = \sum_{i=0}^ka_i x^i , \label{6.8} \end{equation}
and the associated quadratic form  
\begin{equation} \begin{split} Q=Q_{\mathbf{a}}(\lambda,\rho) &= \langle P_{\mathbf{a}}(x),P_{\mathbf{a}}(x)\rangle_\rho \\& =\sum_{0\le i,j\le k}a_ia_j\langle x -\rho,x^{i+j}\rangle \\&
= \sum_{0\le i,j\le k}a_ia_j\left( \mu_{i+j+1}(\lambda)-\rho \mu_{i+j}(\lambda)\right)  \\&
= \sum_{0\le i,j\le k}a_ia_jc_{i+j}, \end{split}\label{6.9}\end{equation}
where we define
\begin{equation} c_m=c_m(\lambda,\rho)= \mu_{m+1}(\lambda) -\rho \mu_m(\lambda) .\label{6.10}\end{equation}
The optimization problem we need to solve is to maximize $Q$ over all vectors normalized by $a_k=1$ when $\rho >0$ is fixed. The solution involves
the (generalized) Laguerre polynomials defined for $\alpha > -1$ by
\begin{equation} {L_n}^{(\alpha)}(x) = \sum_{\nu=0}^n(-1)^\nu\genfrac{(}{)}{0pt }{0}{n+\alpha}{n-\nu} \frac{x^\nu}{\nu !} .\label{6.11}\end{equation}
The zeros of the Laguerre polynomials are real, positive, and simple, (see Chapter 6 of \cite{Sz}.) We denote the smallest zero of  ${L_n}^{(\alpha)}(x)$ by $x_1(n,\alpha)$.  The solution of our problem is obtained in the following proposition.
\begin{proposition} For each $k\ge 1$ and  $\rho >k $ fixed, we have for $0<\lambda< x_1(k+1,\rho-k-1)$ 
\begin{equation} \max_{a_k=1} Q_{\mathbf{a}}(\lambda,\rho)= - (k+1)! \lambda^k\frac{{L_{k+1}}^{(\rho-k-1)}(\lambda)}{{L_{k}}^{(\rho-k)}(\lambda)}.\label{6.12} \end{equation}
Thus, for each $k\ge 1$ and $\rho >k$
\begin{equation} \inf\left\{ \lambda>0 : Q_{\mathbf{a}}(\lambda,\rho)>0, a_k=1\right\}= x_1(k+1,\rho-k-1). \label{6.13}\end{equation}
\end{proposition}

The proof of this proposition will ultimately reduce to evaluating the determinant
\begin{equation} 
D_k = \det \begin{vmatrix} c_0& c_1 &c_2 &\dots & c_k \\
c_1& c_2 &c_3 &\dots & c_{k+1} \\
c_2& c_3 &c_4 &\dots & c_{k+2} \\
\vdots & \vdots & \vdots & \vdots &\vdots \\
c_k& c_{k+1} &c_{k+2} &\dots & c_{2k} 
\end{vmatrix} \ = \det \left[ c_{i+j}\right]_{\substack{i=0,1,2, \ldots ,k\\ j=0,1,2,\ldots ,k}}.
\label{6.14} \end{equation}
The solution of the optimization problem can be obtained by choosing $\mathbf{a}$ so that $P_{\mathbf{a}}(x)$ is orthogonal to all lower degree polynomials with respect to  $\langle\ ,\ \rangle_\rho$. Thus we consider the $k$ equations
\begin{equation} \langle P_{\mathbf{a}}(x), x^i\rangle_\rho=0, \quad i=0,1,2,\ldots , k-1 , \label{6.15}\end{equation}
and prove the following lemma. 
\begin{lemma} If  $D_{k-1}\neq 0$ for a given $\lambda$, then there is an (explicitly obtained) vector $\mathbf{a}$ with $a_k=1$ which satisfies \eqref{6.15} and for which
\begin{equation}   Q_{\mathbf{a}}(\lambda,\rho)= \frac{\ D_k}{D_{k-1}}. \label{6.16} \end{equation} 
\end{lemma}

\noindent \emph{Proof.} We take $a_k=1$. Equation \eqref{6.15} is equivalent to the equations
\begin{equation} \sum_{j=0}^ka_j c_{i+j}=0,\quad i=0,1,2,\ldots , k-1. \label{6.17} \end{equation}
If $\mathbf{a}$ satisfies these equations, then with $\delta_{ij}$ denoting the Kronecker delta, we have 
\begin{equation} \begin{split} Q &= \sum_{i=0}^k a_i \left(\sum_{j=0}^ka_jc_{i+j}\right)
\\ &= \sum_{i=0}^k a_i \left(\delta_{ik}\sum_{j=0}^ka_jc_{i+j}\right) \\&
=\sum_{j=0}^ka_j c_{j+k} .\label{6.18} \end{split}\end{equation}
On rewriting \eqref{6.17} in the form
\begin{equation} \begin{split}
&c_0a_0 +c_1a_1+c_2a_2+ \cdots + c_{k-1}a_{k-1} = -c_k \\ &
c_1a_0 +c_2a_1+c_3a_2+ \cdots + c_{k}a_{k-1} = -c_{k+1} \\ &
\cdots \cdots \cdots \cdots \cdots \cdots \cdots \cdots \cdots \cdots \cdots \cdots \cdots \cdots \cdots \cdots \\&
c_{k-1}a_0 +c_{k}a_1+c_{k+1}a_2+ \cdots + c_{2k-2}a_{k-1} = -c_{2k-1} \label{6.19} \end{split} \end{equation}
we have by Cramer's rule (see \cite{Us}) that these equations have the solution
\begin{equation} a_j = - \frac{D_{k-1}^{(j+1)}}{D_{k-1}}, \quad j=0,1,\ldots , k-1, \label{6.20}\end{equation} 
provided that $D_{k-1}\neq 0$, where $D_{k-1}^{(i)}$ is the determinant with the $i$th column of $D_{k-1}$ replaced by the column $(c_k,c_{k+1},\ldots , c_{2k-1})$.
Thus \eqref{6.18} gives with this choice 
\begin{equation} Q = \frac{1}{D_{k-1}}\left(-\sum_{j=0}^{k-1}D_{k-1}^{(j+1)}c_{k+j}+D_{k-1}c_{2k}\right) .\label{6.21}\end{equation}
On the other hand,  if we expand $D_k$ into its cofactor expansion along the bottom row we see
\[ D_k= \sum_{j=0}^{k}(-1)^{k+j}D_{k+1,j+1}c_{k+j}, \]
where the minor $D_{i,j}$ is the determinant of the matrix where the $i$th row and $j$th column of $D_k$ is removed. From \eqref{6.14} we see 
\[ D_{k+1,j+1}= (-1)^{k-j-1}D_{k-1}^{(j+1)}, \quad 1\le j\le k-1,\quad  D_{k+1, k+1}=D_{k-1} \]
where the factor $(-1)^{k-j-1}$ results from shifting the last column of $D_k$ by $k-j-1$ places to the left. Hence we conclude
\begin{equation} Q = \frac{D_k}{D_{k-1}}. \label{6.22}\end{equation}

Our next lemma evaluates $D_k$. 
\begin{lemma} We have 
\begin{equation}  D_{k-1} = (-1)^k 1!\, 2!\, 3!\, \cdots \, k!\,\lambda^{\frac{k(k-1)}{2} }{L_k}^{(\rho -k)}(\lambda) . \label{6.23}\end{equation}
\end{lemma}

\emph{Proof.} We first claim that
\begin{equation} D_{k-1} = (-1)^k E_k \label{6.24}\end{equation}
where
\begin{equation}
E_k = \det \begin{vmatrix} \mu_0& \mu_1 &\mu_2 &\dots & \mu_k \\
\mu_1& \mu_2 &\mu_3 &\dots & \mu_{k+1} \\
\mu_2& \mu_3 &\mu_4 &\dots & \mu_{k+2} \\
\vdots & \vdots & \vdots & \vdots &\vdots \\
\mu_{k-1}& \mu_{k} &\mu_{k+1} &\dots & \mu_{2k-1} \\
1&\rho&\rho^2 &\ldots &\rho^k 
\end{vmatrix} \ = \det \begin{vmatrix} \mu_{i+j}\\ \rho^j\end{vmatrix}_{\substack{i=0,1,2, \ldots ,k-1\\ j=0,1,2,\ldots ,k}},
\label{6.25}\end{equation}
for if in $E_k$ we multiply the $\ell$-th column by $\rho$ and subtract this from the $(\ell+1)$-th column for $\ell =1,2,\ldots , k$ we obtain 
 \[
E_k = \det \begin{vmatrix} \mu_0& c_0 &c_1 &\dots & c_{k-1} \\
\mu_1& c_1 &c_2 &\dots & c_{k} \\
\mu_2& c_2 &c_3 &\dots & c_{k+1} \\
\vdots & \vdots & \vdots & \vdots &\vdots \\
\mu_{k-1}& c_{k-1} &c_{k} &\dots & c_{2k-2} \\
1&0&0 &\ldots &0 
\end{vmatrix},\]
and using the cofactor expansion along the bottom row gives $E_k= (-1)^k D_{k-1}$. 

We now introduce the differential operators
\begin{equation} D = \frac{d}{d\lambda}, \quad  \delta = \lambda D = \lambda \frac{d}{d\lambda}, \quad \Delta = \delta + \lambda =\lambda \frac{d}{d\lambda}+\lambda .
\label{6.26}\end{equation}
Clearly we have the relations
\begin{equation} \delta^k =  \lambda^kD^k+\sum_{j=1}^{k-1}a_j(\lambda) D^j, \quad \Delta^k = \delta^k + \sum_{j=0}^{k-1}b_j(\lambda)\delta^j \label{6.27}\end{equation}
where $a_j(\lambda)$ and $b_j(\lambda)$ are polynomials of degree $j$ in $\lambda$. Now by \eqref{6.4} and \eqref{6.5} we have $\mu_k = \Delta \mu_{k-1}$, and in general
\begin{equation} \mu_k = \Delta^i \mu_{k-i}\ ,\  0\le i\le k; \qquad  \mu_k = \Delta^k 1 .\label{6.28} \end{equation}
From this we see
\begin{equation} E_k = \det \begin{vmatrix} \mu_0& \mu_1 &\mu_2 &\dots & \mu_k \\
\Delta \mu_0& \Delta \mu_1 &\Delta \mu_2 &\dots & \Delta\mu_{k}\\
\Delta^2\mu_0& \Delta^2\mu_1 &\Delta^2\mu_3 &\dots &\Delta^2 \mu_{k} \\
\vdots & \vdots & \vdots & \vdots &\vdots \\
\Delta^{k-1}\mu_{0} &\Delta^{k-1}\mu_{1} &\Delta^{k-1}\mu_{2} &\dots & \Delta^{k-1}\mu_{k} \\ 1& \rho & \rho^2 & \dots & \rho^k
\end{vmatrix}  \ = \det \begin{vmatrix}\Delta^i \mu_{j}\\ \rho^j\end{vmatrix}_{\substack{i=0,1,2, \ldots ,k-1\\ j=0,1,2,\ldots ,k}}.\label{6.29}\end{equation}
By the second relation in \eqref{6.27} we can replace $\Delta^i$ by $\delta^i$ and a linear combination of lower powers of $\delta$, which can be elimated by row operations. Thus we can replace $\Delta$ by $\delta$ in the above determinant without effecting its value, and then by the first relation in \eqref{6.27} and row operations we can replace $\delta^i$ by $\lambda^i D^i$ which on removing the factors of $\lambda$ in each row gives
 \begin{equation} \begin{split} E_k &= \lambda^{\frac{k(k-1)}{2}}\det \begin{vmatrix} \mu_0& \mu_1 &\mu_2 &\dots & \mu_k \\
D\mu_0& D\mu_1 &D\mu_2 &\dots & D\mu_{k}\\
D^2\mu_0& D^2\mu_1 &D^2\mu_3 &\dots & D^2\mu_{k} \\
\vdots & \vdots & \vdots & \vdots &\vdots \\
D^{k-1}\mu_{0} &D^{k-1}\mu_{1} &D^{k-1}\mu_{2} &\dots &D^{k-1} \mu_{k} \\ 1& \rho & \rho^2 & \dots & \rho^k
\end{vmatrix} \\ &= \lambda^{\frac{k(k-1)}{2}}\det \begin{vmatrix}D^i \mu_{j}\\ \rho^j\end{vmatrix}_{\substack{i=0,1,2, \ldots ,k-1\\ j=0,1,2,\ldots ,k}}.\end{split}
\label{6.30}\end{equation}
We next need the relation
\begin{equation} q(q-1)\cdots (q-h+1)= h! \left(\begin{array}{c} q\\ h  	
\end{array} \right) = \sum_{j=0}^h(-1)^{h-j}\left[ \begin{array}{c} h\\ j  	
\end{array} \right]q^j ,\label{6.31}\end{equation}
where $\left[ \begin{array}{c} h\\ j  	
\end{array} \right]$ are the Stirling numbers of the first type, although we do not need to use any properties of these numbers. Then we have by \eqref{6.3}
\begin{equation} \begin{split} \lambda^h = e^\lambda \lambda^he^{-\lambda} &= \sum_{q=0}^\infty \frac{\lambda^{q+h}}{q!}e^{-\lambda} = \sum_{q=0}^\infty q(q-1)\cdots (q-h+1)\frac{\lambda^q}{q!}e^{-\lambda}  \\ &= \sum_{q=0}^\infty \left( \sum_{j=0}^h(-1)^{h-j}\left[ \begin{array}{c} h\\ j  	
\end{array} \right]q^j \right)\frac{\lambda^q}{q!}e^{-\lambda}= \sum_{j=0}^h (-1)^{h-j}\left[ \begin{array}{c} h\\ j  	
\end{array} \right]\mu_j(\lambda) .\end{split}\label{6.32} \end{equation}
Thus, using column operations we see
\begin{equation} \begin{split} E_k &= \lambda^{\frac{k(k-1)}{2}}\det \begin{vmatrix} 1& \lambda &\lambda^2 &\dots & \lambda^k \\
D 1& D\lambda &D\lambda^2 &\dots & D\lambda^{k}\\
D^21& D^2\lambda &D^2\lambda^2 &\dots & D^2\lambda^k \\
\vdots & \vdots & \vdots & \vdots &\vdots \\
D^{k-1}1 &D^{k-1}\lambda &D^{k-1}\lambda^2 &\dots &D^{k-1} \lambda^k \\ 1& 1! \left(\begin{array}{c} \rho \\ 1 	
\end{array} \right) & 2! \left(\begin{array}{c} \rho\\ 2  	
\end{array} \right)& \dots & k! \left(\begin{array}{c} \rho\\ k  	
\end{array} \right)
\end{vmatrix} \\ &= \lambda^{\frac{k(k-1)}{2}}\det \begin{vmatrix}D^i \lambda^j\\ j! \left(\begin{array}{c} \rho\\ j 	
\end{array} \right)\end{vmatrix}_{\substack{i=0,1,2, \ldots ,k-1\\ j=0,1,2,\ldots ,k}}.\end{split}
\label{6.33}\end{equation}
Expanding along the bottom row we see
\begin{equation} \det \begin{vmatrix}D^i \lambda^j\\ j! \left(\begin{array}{c} \rho\\ j 	
\end{array} \right)\end{vmatrix}_{\substack{i=0,1,2, \ldots ,k-1\\ j=0,1,2,\ldots ,k}}= \sum_{h=0}^k (-1)^{k-h}h!\left(\begin{array}{c} \rho\\ h	
\end{array} \right)\det[D^i\lambda^j]_{\substack{i=0,1,\ldots,k-1 \\ j=0,1,\ldots , k; j\neq h}}. \label{6.34}\end{equation}
We will show below that 
\begin{equation}F_h = \det[D^i\lambda^j]_{\substack{i=0,1,\ldots,k-1 \\ j=0,1,\ldots , k; j\neq h}}=1!\, 2!\, \cdots \, (k-1)! \left(\begin{array}{c} k\\ h	
\end{array} \right) \lambda^{k-h}\label{6.35}\end{equation}
which then gives on retracing our steps
\[\begin{split} D_{k-1} &=(-1)^k \lambda^{\frac{k(k-1)}{2}}1!\, 2!\, \cdots\,  (k-1)!\sum_{h=0}^k (-1)^{k-h}h!\left(\begin{array}{c} \rho\\ h	
\end{array} \right) \left(\begin{array}{c} k\\ h	
\end{array} \right) \lambda^{k-h} \\& = (-1)^k \lambda^{\frac{k(k-1)}{2}}1!\, 2!\, \cdots \, k! {L_k}^{(\rho-k)}(\lambda),\end{split}\]
which proves Lemma 4. 

We prove \eqref{6.35} by the following argument shown to us by Wasin So. We consider the complete upper triangular matrix
\[ \begin{split} M &= [D^i\lambda^j]_{\substack{i=0,1,\ldots,k\\ j=0,1,\ldots ,k}} = \begin{pmatrix} 1& \lambda &\lambda^2 &\dots & \lambda^k \\
D 1& D\lambda &D\lambda^2 &\dots & D\lambda^{k}\\
D^21& D^2\lambda &D^2\lambda^2 &\dots & D^2\lambda^k \\
\vdots & \vdots & \vdots & \vdots &\vdots \\
D^{k}1 &D^{k}\lambda &D^{k}\lambda^2 &\dots &D^{k} \lambda^k
\end{pmatrix}\\ & = \left[i! \binom{j}{i}\lambda^{j-i}\right]_{\substack{i=0,1,\ldots,k\\ j=0,1,\ldots, k}}.\end{split}
\]
Observe that $\det M = 1!\, 2!\, \cdots k!$, and further
that 
\[ M = TP,\quad  \text{where }\ T= [ \delta_{ij} i!\ ]_{\substack{i=0,1,\ldots,k\\ j=0,1,\ldots k}}, \quad P = \left[ \binom{j}{i}\lambda^{j-i}\right]_{\substack{i=0,1,\ldots,k\\ j=0,1,\ldots ,k}}.\]
Now by the matrix inverse formula using minors 
\[ M^{-1} = \frac{1}{\det M}[ (-1)^{i+j}D_{j,i}]_{\substack{i=0,1,\ldots,k\\ j=0,1,\ldots ,k}},\]
where $F_h$ occurs in this matrix as the minor $D_{k,h}$. Further,
\[ M^{-1} = P^{-1}T^{-1},\]
where
\[ T^{-1}= [ \delta_{ij}\frac{1}{ i!}\ ]_{\substack{i=0,1,\ldots,k\\ j=0,1,\ldots ,k}},\quad \ P^{-1} = \left[ (-1)^{j-i}\binom{j}{i}\lambda^{j-i}\right]_{\substack{i=0,1,\ldots,k\\ j=0,1,\ldots ,k}}, \]
where we used the identity
\[ \sum_{s=0}^k (-1)^s\binom{j}{s}\binom{s}{i} = (-1)^i\delta_{ij} .\]
From this last relation we see, letting $M^{-1}=[\bar{m}_{ij}]$, 
\[ F_h = (-1)^{k+h}(\det M) \bar{m}_{hk} = (-1)^{k+h}\frac{(-1)^{k-h}}{k!} 1!\, 2!\, \cdots\,  k! \binom{k}{h}\lambda^{k-h},\]
as desired. 

\emph{Proof of Proposition 3}. Let $\mathbf{a}$ be the solution for \eqref{6.15} found in Lemma 3, which exists for any  $\lambda$ where $D_{k-1}\neq 0$, and let $\mathbf{b}$ be any other $k$-vector with $b_k=1$.  Then   $P_{ \mathbf{b-a}}(x)$ is a polynomial of degree $k-1$ or less, and by the orthogonality property \eqref{6.15}
\[ \begin{split} Q_{\mathbf{b}}(\lambda,\rho) &= \langle P_{\mathbf{a}}(x)+P_{\mathbf{b-a}}(x),P_{\mathbf{a}}(x)+P_{\mathbf{b-a}}(x)\rangle_\rho \\&
=\langle P_{\mathbf{a}}(x),P_{\mathbf{a}}(x)\rangle_\rho +\langle P_{\mathbf{b-a}}(x),P_{\mathbf{b-a}}(x)\rangle_\rho\\&  =Q_{\mathbf{a}}(\lambda,\rho)+ Q_{\mathbf{b-a}}(\lambda,\rho).
\end{split} \]
In general by \eqref{6.7} and \eqref{6.9} for any $\mathbf{c} \neq \mathbf{0}$, assuming $\rho>0$ is fixed,
\[ \begin{split} Q_{\mathbf{c}}(\lambda,\rho)&= \sum_{j=0}^\infty (j-\rho)(P_{\mathbf{c}}(j))^2p(j)\\&
= -\rho \langle (P_{\mathbf{c}}(x))^2 , 1\rangle + O_{\mathbf{c}}(\lambda) \\&
< 0, \quad \text{for} \quad 0<\lambda \le \lambda_0(\mathbf{c},\rho)\end{split}\]
where $\lambda_0(\mathbf{c},\rho)$ is a small positive constant depending on $\mathbf{c}$.  Thus
\[ Q_{\mathbf{b}}(\lambda,\rho)  \le Q_{\mathbf{a}}(\lambda,\rho)\]
for $0<\lambda <\lambda_0(\mathbf{c})$, 
proving that $Q_{\mathbf{a}}$ is  maximal at least for small enough $\lambda$. This will continue to be true for larger $\lambda$ as long as $Q_{\mathbf{c}}<0$ for any $(k-1)$-vector $\mathbf{c}$, and therefore as long as the maximal $Q$ for $(k-1)$-vectors is negative. By
(5.1.14) of Szeg\"o \cite{Sz} we have
\[ \frac{d}{dx}{L_n}^{(\alpha)}(x)=-{L_{n-1}}^{(\alpha+1)}(x),\]
and therefore we see the sequence $\{{L_k}^{(\rho-k)}\}$ of Laguerre polynomials has the property that the negative of the derivative of a term is the previous term. (Thus the negative derivative of the Laguerre polynomial in the numerator in \eqref{6.12} is the Laguerre polynomial in the denominator.) Further this sequence of Laguerre polynomials all are decreasing functions up to their first positive zero, and hence the sequence of smallest positive zeros $x_1(k,\rho-k)$ is a decreasing sequence. Starting with the trivial case when $k=1$ we see successively that the  $Q_k$  with $\mathbf{a}$ satisfying \eqref{6.15} will be maximal for 
$0<\lambda< x_1(k+1,\rho-k-1)$. This completes the proof of Proposition 3. 

Our next result evaluates the smallest positive zero $x_1(n,\alpha)$ asymptotically as $n\to \infty$.  
\begin{lemma}  Let ${L_n}^{(\alpha)}(x)$, $\alpha > -1$, denote the Laguerre polynomials. The zeros of ${L_n}^{(\alpha)}(x)$ are real, positive, and simple. Let $x_1(n,\alpha)$ denote the smallest zero of ${L_n}^{(\alpha)}(x)$.   If $\alpha = \beta(n)-n$ and  $\lim_{n\to \infty}\frac{ \beta(n) }{n} = A>0$, then 
\begin{equation}  \lim_{n\to\infty}\frac{x_1(n,\alpha)}{n} = (\sqrt{A}-1)^2. \label{3.15}\end{equation}
\end{lemma} 
\bigskip
\noindent \emph{Proof.} The properties of ${L_n}^{(\alpha)}(x)$ may be found in Szeg\"o \cite{Sz}. Equation \eqref{3.15} is a special case of Theorem 4.4 of \cite{DS}. A simple proof may be obtained by using the same argument found in \cite{MSV} where a result corresponding to \eqref{3.15} for Jacobi polynomials is proved using Sturm comparison theory.  By (5.1.2) of Szeg\"o, the differential equation
\begin{equation} u'' +\left(\frac{n+(\alpha+1)/2}{x} +\frac{1-\alpha^2}{4x^2}- \frac{1}{4}\right) u =0 \label{3.16} \end{equation}
has $u=e^{-x/2}x^{(\alpha+1)/2}{L_n}^{(\alpha)}(x)$ as a solution.
Let
\[\begin{split} H_n(x) :&=  \frac{n+(\alpha+1)/2}{x} +\frac{1-\alpha^2}{4x^2}- \frac{1}{4}\\ & = - \frac{ x^2 -(4n+2(\alpha+1))x +(\alpha^2-1)}{4x^2}, \end{split} \]
and denote the smaller root of the quadratic in the numerator by ${x_n}^-$. Then by the Sturm comparison argument in \cite{MSV}, and noting $\lim_{n\to\infty}\frac{\alpha}{n}  = A-1$, we have
\[ \begin{split} \lim_{n\to\infty}\frac{x_1(n,\alpha)}{n} &= \lim_{n\to \infty} \frac{{x_n}^-}{n} \\ &= \lim_{n\to \infty}\frac{ 2n +\alpha +1 -\sqrt{4n^2 +2\alpha +2 +4n\alpha +4n}}{n} \\ &= 2 + \lim_{n\to\infty}\frac{\alpha}{n}- \sqrt{ 4 + 4\lim_{n\to\infty}\frac{\alpha}{n}} \\ & = (\sqrt{A}-1)^2.\end{split}\]
\section{Gaps between primes.}
In this section we prove Theorem 1. We want to examine statistically the number of primes in the interval $(n,n+h]$ for $N <n \le 2N$ with $N\to \infty$. In this range the average distance between consecutive primes is $\log N$, and thus we will take $h$ to be a multiple of this length. We therefore let
\begin{equation} \psi(x) = \sum_{n\le x}\Lambda(n),  \label{7.1}\end{equation}
\begin{equation} \psi(n,h) = \psi(n+h)-\psi(n) ,  \label{7.2} \end{equation}
\begin{equation}   h=\lambda \log N , \label{7.3}\end{equation}
and in this paper we assume that
\begin{equation}\lambda \ll 1 . \label{7.4} \end{equation}
The model for our method is due to Gallagher \cite{GA}, who proved that if the Hardy-Littlewood conjecture \eqref{2.4} holds uniformly for  $ h \ll \log N$ then one can asymptotically evaluate all the moments for the number of primes in intervals of length $h$. Thus assuming \eqref{2.4}, Gallagher proved that  
\begin{equation}  M_k(N,h,\psi):= \frac{1}{N(\log N)^k}\sum_{n=N+1}^{2N} (\psi(n,h))^k \sim \mu_k(\lambda),\label{7.5} \end{equation}
as $N\to \infty$, where $\mu_k(\lambda)$ is 
the Poisson moment from \eqref{6.3} and \eqref{6.4}. 

In order to obtain unconditional results we make use of our approximation $\Lambda_R(n;\mathcal{H})$. Taking $N<n\le 2N$, we first need to approximate 
\begin{equation} \begin{split} \psi(n,h)^k &= \sum_{1\le h_1,h_2,\ldots ,h_k\le h}\Lambda(n+h_1)\Lambda(n+h_2)\cdots \Lambda(n+h_k)\\& = (1+o(1))\sum_{1\le h_1,h_2,\ldots ,h_k\le h}(\log N)^{k-|\mathcal{H}|}\Lambda(n;\mathcal{H}). \end{split}\label{7.6}\end{equation}
To define our approximation, we extend the definition of $\Lambda_R(n;\mathcal{H})$ in \eqref{2.3} to vectors (or lists)
$\mbf{H}=(h_1,h_2,\ldots, h_k)$. The distinct components of the vector   $\mbf{H}$ are the elements of the set $\mathcal{H}=\{h_1,h_2,\ldots , h_k\}$, and we define 
\begin{equation} \Lambda_R(n;\mbf{H}) := (\log R)^{k-|\mathcal{H}|}\Lambda_R(n;\mathcal{H}).\label{7.7} \end{equation}
Then our approximation of $\psi(n,h)^k$ is
\begin{equation} {\psi_R}^{(k)}(n,h) := \sum_{1\le h_1,h_2,\ldots , h_k\le h}\Lambda_R(n;\mbf{H}).\label{7.8}\end{equation}
For convenience we also define ${\psi_R}^{(0)}(n,h)=1$.
We next define the approximate moments, letting $k=i+j$,
\begin{equation}  M_{ij}(R) = \frac{1}{N(\log R)^k}\sum_{n=N+1}^{2N} {\psi_R}^{(i)}(n,h){\psi_R}^{(j)}(n,h) ,\label{7.9} \end{equation}
and note that $M_{00}(R)= 1$. We also need the mixed moments
\begin{equation} \tilde{ M}_{ij}(R) = \frac{1}{N(\log R)^{k+1}}\sum_{n=N+1}^{2N} {\psi_R}^{(i)}(n,h){\psi_R}^{(j)}(n,h)\psi(n,h) ,\label{7.10} \end{equation}
for which we note by the prime number theorem that $\tilde{ M}_{00}(R)\sim \frac{\lambda}{\theta} = \mu_1(\frac{\lambda}{\theta})$, in accord with \eqref{7.11} and \eqref{7.13} below.

Using Propositions 1 and 2 we will prove asymptotically that these approximate moments are also Poisson moments with an increased expected value involving the truncation level $R$. Define $\theta$ by
\begin{equation} R=N^\theta .\label{7.11}\end{equation}
\begin{proposition}  As $N\to \infty$ we have, for $k=i+j \ge 1$ and  for any fixed $0<\theta < \frac{1}{k}$, 
\begin{equation} M_{ij}(R) =(1+o_{k}(1)) \mu_{k}\left(\frac{\lambda}{\theta}\right) , \label{7.12}\end{equation}
and for any fixed $0<\theta < \frac{1}{2k}$,
\begin{equation} \tilde{M}_{ij}(R) = (1+o_{k}(1))\mu_{k+1}\left(\frac{\lambda}{\theta}\right) . \label{7.13}\end{equation}
\end{proposition}
\emph{Proof.}  By differencing, Propositions 1 and 2 continue to hold unchanged when we sum for $N<n\le 2N$. We first extend Proposition 1 for vectors $\mbf{H}_1$ and $\mbf{H}_2$. Recalling the notation $|\mbf{H}|$ which denotes the number of components of the vector $\mbf{H}$, let $k= |\mbf{H}_1|+ |\mbf{H}_2|$, $\mathcal{H}= \mathcal{H}_1\cup \mathcal{H}_2$, where $\mathcal{H}_i$ is the set of distinct components of $\mbf{H}_i$.  Then by \eqref{7.7} and Proposition 1, (note the $k$ in Proposition 1 and 2 is equal to $|\mathcal{H}_1| +|\mathcal{H}_2|$ here), we have for $R=o(N^{\frac{1}{k}})$,  
\begin{equation}\begin{split} \sum_{n=N+1}^{2N}\Lambda_R(n;\mbf{H}_1)&\Lambda_R(n;\mbf{H}_2)\\& = (\log R)^{k-|\mathcal{H}_1|-|\mathcal{H}_2|}\sum_{n=N+1}^{2N}\Lambda_R(n;\mathcal{H}_1)\Lambda_R(n;\mathcal{H}_2)\\& = N\big(\gs(\mathcal{H})+o(1)\big)(\log R)^{k-|\mathcal{H}|}.\label{7.14}\end{split} \end{equation}
Thus we see this result depends on $k$ and not the individual values of $|\mbf{H}_1|$ and $|\mbf{H}_2|$.
Hence, letting $h_1,h_2,\ldots, h_k$ list the components of $\mbf{H}_1$ and $\mbf{H}_2$ (in any order), we have
\[ \begin{split} M_{ij}(R) &= \frac{1}{N(\log R)^k}\sum_{1\le h_1,h_2,\ldots , h_k\le h} \sum_{n=N+1}^{2N}\Lambda_R(n;\mbf{H}_1)\Lambda_R(n;\mbf{H}_2) \\& = \sum_{1\le h_1,h_2,\ldots , h_k\le h} \big(\gs(\mathcal{H})+o_k(1)\big)(\log R)^{-|\mathcal{H}|},
\end{split}\]
provided $R=o(N^{\frac{1}{k}})$.
We group terms in this sum according to the number of distinct values $\nu$ of $h_1,h_2,\ldots, h_k$, and denote these distinct values by ${h'}_1, {h'}_2, \ldots, {h'}_\nu$.  There are $\genfrac{\{}{\}}{0pt }{1}{k}{\nu}$ ways to partition the $k$ $h_i$'s into these $\nu$ disjoint sets, and all of these will occur in the sum above.
Hence by \eqref{2.13} we have
\[ \begin{split} M_{ij}(R) &=  \sum_{\nu=1}^k \genfrac{\{}{\}}{0pt }{0}{k}{\nu} h^\nu (1+o_k(1))(\log R)^{-\nu}  \\
& =(1+o_k(1)) \mu_k\left(\frac{\lambda}{\theta}\right) \end{split}\]
which proves the first part of Proposition 4. The second part is proved identically using Proposition 2.

Now consider  
\begin{equation}\begin{split}\mathcal{S}_{k}&= \mathcal{S}_{k}(N,R,\lambda , \rho) \\&= \frac{1}{N(\log R)^{2k+1}}\sum_{n=N+1}^{2N} (\psi(n,h) - \rho \log N)\big(P_k(\psi_R(n,h) )\big)^2 ,\end{split}\label{7.15}\end{equation}
where
\begin{equation}  P_k(\psi_R(n,h) ) = \sum_{\ell=0}^k a_\ell {\psi_R}^{(\ell)}(n,h) (\log R)^{k-\ell},\label{7.16}\end{equation}
and the $a_\ell$'s are arbitrary functions of $N$, $R$, $k$, $\lambda$, and $\rho$ which are to be chosen to optimize the argument. 
On multiplying out we have that
\begin{equation}\begin{split}\mathcal{S}_{k} &= \frac{1}{N(\log R)^{2k+1}}\sum_{0\le i,j\le k} a_ia_j(\log R)^{2k-i-j}\\&\hskip 1.3in \times \sum_{n=N+1}^{2N} (\psi(n,h) - \frac{\rho}{\theta} \log R){\psi_R}^{(i)}(n,h){\psi_R}^{(j)}(n,h) \\&
= \sum_{0\le i,j\le k} a_ia_j\mathcal{M}_{ij}. 
\end{split}\label{7.17}\end{equation}
Letting  
\begin{equation}\tilde{\lambda}= \frac{\lambda}{\theta}, \qquad \tilde{\rho}=\frac{\rho}{\theta},\label{7.18}\end{equation} 
we have by Proposition 4 on taking $i+j =\kappa$ and assuming $0<\theta < \frac{1}{2\kappa}$,
\begin{equation} \begin{split}\mathcal{M}_{ij} &= \tilde{M}_{ij}(R)- \frac{\rho}{\theta} M_{ij}(R)\\& = \mu_{\kappa+1}(\tilde{\lambda}) - \tilde{\rho} \mu_{\kappa}(\tilde{\lambda})+ o_\kappa(\frac{1}{ \theta^{\kappa+1}}) \\& = c_{\kappa}(\tilde{\lambda},\tilde{\rho})+ o_\kappa (\frac{1}{ \theta^{\kappa+1}}),\label{7.19}\end{split}\end{equation}
using the notation of \eqref{6.10} in the last line. 
To evaluate $\mathcal{S}_k$ we need to apply these results for $0\le \kappa \le 2k$, all of which will hold if we impose the condition 
\begin{equation}  \frac{1}{4k +1} \le \theta <\frac{1}{4k}.\label{7.20}\end{equation} 
Thus
\begin{equation} \begin{split}\mathcal{S}_k &= \sum_{0\le i,j\le k} a_ia_jc_{i+j}(\tilde{\lambda},\tilde{\rho}) +o_k(\max_{1\le \ell \le k}|a_\ell|^2)\\& = Q_{\mathbf{a}} (\tilde{\lambda},\tilde{\rho}) +o_k(1), \label{7.21}\end{split} \end{equation}
since  $\max_{1\le \ell \le k}|a_\ell|^2$ depend only on $k$ for fixed $\lambda$ and $\rho$.
By Proposition 3 we obtain a sign change for $Q_{\mathbf{a}} (\tilde{\lambda},\tilde{\rho}) $ at the smallest zero  $x_1(k+1,\tilde{\rho}-k-1)$ of the Laguerre polynomial ${L_{k+1}}^{(\tilde{\rho}-k-1)}(\tilde{\lambda})$, with $Q_{\mathbf{a}} (\tilde{\lambda},\tilde{\rho}) $  negative  for $0< \tilde{\lambda} < x_{1}(k+1,\tilde{\rho}-k-1)$ and positive for $ x_{1}(k+1,\tilde{\rho}-k-1)< \tilde{\lambda} < x_{1}(k,\tilde{\rho}-k)$.  Therefore by \eqref{7.21}  $\mathcal{S}_k$ will also be positive for $ x_{1}(k+1,\tilde{\rho}-k-1)+o_k(1)< \tilde{\lambda} < x_{1}(k,\tilde{\rho}-k)- o_k(1)$ as $N\to \infty$.  We apply Lemma 5 with $\beta(k) = \tilde{\rho}$; if we take sequences $\theta =\theta_k \to {\frac{1}{4k}}^-$ and $\rho = \rho_k \to r^+$ as  $k\to \infty$, then $A= 4r$,
 and there exists constants $0<c_k <{c'}_k$, and $c_k,{c'}_k \to 0$, such that for  $(\sqrt{r}-\frac{1}{2})^2+c_k \le \lambda \le (\sqrt{r}-\frac{1}{2})^2+ {c'}_k$ we have
\begin{equation} \mathcal{S}_k \gg_k 1, \quad   \mathcal{S}_k>0.\label{7.22}\end{equation}
Note that the Laguerre polynomials are well defined by \eqref{7.20} here since by \eqref{7.20}  $\tilde{\rho}-k>0$.
The proof of Theorem 1 is now a standard deduction from \eqref{7.22}; we follow our earlier proof in the last section of \cite{GYI}. Define
\begin{equation} Q_r^+(N,h) = \sum_{\substack{n=N+1\\ \pi(n+h)-\pi(n)> r}}^{2N} 1.
\label{7.23}
\end{equation}
 If $n$ is an integer for which $\pi(n+h)-\pi(n) >r$ then there must be a $j$ such that $n\le p_j$ and $p_{j+r}\le n+h$.
Thus $p_{j+r}-p_{j}\le h$ and $p_{j+r}-h\le n \le p_j<p_{j+r}$, so that there are at most $h$ such $n$'s corresponding to each such gap. Therefore
\begin{equation} Q_r^+(N,h) \ll_r h \sum_{\substack{ N<p_n \le 2N\\ p_{n+r}-p_n \le h}} 1+O(Ne^{-c\sqrt{\log N}}), \label{7.24}\end{equation}
where we have used the prime number theorem to remove the prime gaps overlapping the endpoints $N$ and $2N$. (This can be done more explicitly as in \cite{GYI}.)

Next, we have, for $N$ sufficiently large,
\begin{equation}
Q_r^+(N,h) =  \sum_{\substack{n=N+1\\ \psi(n+h)-\psi(n) \ge \rho \log N}}^{2N} 1 +O(N^{\frac{1}{2}}),\label{7.25}
\end{equation}
where $\rho$ can be taken to be any number in the range $r<\rho < r+1$, and the error term is from removing prime powers.  By \eqref{7.25} and Cauchy's inequality  we see that 
\begin{equation}\begin{split} \mathcal{S}_{k} &\le  \frac{1}{N(\log R)^{2k+1}}\sum_{\substack{n=N+1\\ \psi(n+h)-\psi(n)\ge \rho \log N}}^{2N} \psi(n,h)\big(P_k(\psi_R(n,h) )\big)^2 \\
&\le  \frac{1}{N(\log R)^{2k+1}} \left(\sum_{\substack{n=N+1\\ \psi(n+h)-\psi(n)\ge \rho \log N}}^{2N} \big(P_k(\psi_R(n,h) )\big)^2 \right)^{\frac{1}{2}}\\& \hskip 2in \times \left(\sum_{n=N+1}^{2N} \psi(n,h)^2\big(P_k(\psi_R(n,h) )\big)^2 \right)^{\frac{1}{2}}\\
&\le  \frac{\sqrt[4]{Q_r^+(N,h)+O(N^{1/2})}}{N(\log R)^{2k+1}} \left(\sum_{n=N+1}^{2N}\psi(n,h)^4 \right)^{\frac{1}{4}}\left(\sum_{n=N+1}^{2N}\big(P_k(\psi_R(n,h) )\big)^4 \right)^{\frac{1}{2}}.\label{7.26}\end{split} \end{equation}
Hence, provided \eqref{7.22} holds we have
\begin{equation}  Q_r^+(N,h)+O(N^{1/2}) \gg_k \frac{ \big(N (\log R)^{2k+1}\big)^4}{\displaystyle \left(\sum_{n=N+1}^{2N}\psi(n,h)^4\right)\left(\sum_{n=N+1}^{2N}\big(P_k(\psi_R(n,h) )\big)^4\right)^2} . \label{7.27} \end{equation}
We will prove below that subject to $h\ll \log N$ from \eqref{7.3} and \eqref{7.4} we have
\begin{equation} \sum_{n=N+1}^{2N}\psi(n,h)^4 \ll N(\log N)^4 \label{7.28} \end{equation}
and 
\begin{equation} \sum_{n=N+1}^{2N}\big(P_k(\psi_R(n,h) )\big)^4\ll N(\log N)^{4k}. \label{7.29} \end{equation}
Therefore we conclude from \eqref{7.24} -- \eqref{7.27} that for $\lambda = (\sqrt{r}-\frac{1}{2})^2+c_k $
\begin{equation} \sum_{\substack{ N<p_n \le 2N\\ p_{n+r}-p_n \le h}} 1 \gg_k \frac{N}{h} \gg_k \pi(N) ,\label{7.30} \end{equation}
where $c_k\to 0^+$ as $k\to \infty$, which proves Theorem 1. 

Before proceeding to the proofs of \eqref{7.28} and \eqref{7.29}, we note that, for $N<n\le 2N$, the trivial estimates  $\psi(n,h)\ll h\log N$ and ${\psi_R}^{(k)}(n,h) \ll_k N^\epsilon$ immediately imply the bounds $\ll N^{1+\epsilon}$ in \eqref{7.28} and \eqref{7.29} from which \eqref{1.10} follows.  To prove \eqref{7.28} we make use of the sieve bound 
\begin{equation} \sum_{n\le N}\Lambda(n; \mathcal{H}_k) \le (2^k k! + \epsilon)\gs(\mathcal{H}_k)N, \label{7.31} \end{equation}
(see \cite{HR} Theorem 5.7 or \cite{Gr} Theorem 4 of 2.3.3). Then by  equation \eqref{7.6} and \eqref{2.13}
\begin{equation}\begin{split} \sum_{n=N+1}^{2N}& \psi(n;h)^4 = \sum_{n=N+1}^{2N}(1+o(1))\sum_{1\le h_1,h_2,h_3,h_4\le h}(\log N)^{4-|\mathcal{H}|}\Lambda(n;\mathcal{H})\\&=(1+o(1))\sum_{\nu =1}^4\genfrac{\{}{\}}{0pt }{0}{4}{\nu}(\log N)^{4-\nu}\sum_{\substack{1\le h_1,\ldots ,h_\nu\le h\\ \text{distinct}}}\Big(\sum_{n=N+1}^{2N} \Lambda(n;\mathcal{H}_\nu)\Big)\\&\le (N+o(N))\sum_{\nu =1}^4\genfrac{\{}{\}}{0pt }{0}{4}{\nu}2^\nu \nu ! (\log N)^{4-\nu}\Big(\sum_{\substack{1\le h_1,\ldots ,h_\nu\le h\\ \text{distinct}}}\gs(\mathcal{H}_\nu)\Big)\\&\le \left(\sum_{\nu =1}^4\genfrac{\{}{\}}{0pt }{0}{4}{\nu} \nu ! (2\lambda)^\nu +\epsilon \right) N\log^4 N,\end{split}\label{7.32} \end{equation}
which proves \eqref{7.28}.

The proof of \eqref{7.29} is based on a generalization of Proposition 1 proved in \cite{GYIII} by the same method used in the proof of Proposition 1 in this paper. For 
$k\ge 1$, and $\mathcal{H} =\{h_1, h_2, \ldots , h_r\}$ with distinct integers $h_i$, and $\mbf{a} = (a_1,a_2, \ldots a_r)$,  $a_i\geq 1$ with $\sum_{i=1}^r a_i = k$, let 
\begin{equation} \mathcal{ S}_k(N,\mathcal{H}, \mbf{a}) =\sum_{n=1}^N \Lambda_R(n+h_1)^{a_1}\Lambda_R(n+h_2)^{a_2}\cdots \Lambda_R(n+h_r)^{a_r}.\label{7.33}\end{equation} 
 Then for $\max_i |h_i|\le R $ and $R\ge 2$   we have 
\begin{equation} \mathcal{S}_k(N,\mathcal{H},\mbf{a}) = \big(\mathcal{ C}_k(\mbf{a})\gs(\mathcal{H})+o_k(1)\big)N(\log R)^{k-r} +O(R^k),\label{7.34}\end{equation}
where the $\mathcal{ C}_k(\mbf{a})$ are constants that are computable rational numbers.
On multiplying out the left-hand side of \eqref{7.29} we obtain a linear combination of $(k+1)^4$ terms of the form 
\[ \mathcal{T}(\ell_1,\ell_2,\ell_3,\ell_4)=(\log R)^{4k-\ell_1-\ell_2-\ell_3-\ell_4} \sum_{n=N+1}^{2N}\prod_{i=1}^4 {\psi_{R}}^{(\ell_i)}(n,h) , \] 
for any $0\le \ell_1, \ell_2, \ell_3,\ell_4 \le k$. Letting $\ell = \ell_1+\ell_2+\ell_3+\ell_4$, then $0\le \ell \le 4k$ and we have 
\[ \mathcal{T}(\ell_1,\ell_2,\ell_3,\ell_4)= (\log R)^{4k-\ell}\sum_{1\le m_1,m_2,\ldots , m_\ell\le h}\sum_{n=N+1}^{2N}\prod_{i=1}^4\Lambda_R(n, \mbf{H}_i),\]
where $m_1,m_2,\ldots, m_\ell$ run through the components of the $\mbf{H}_i$, $1\le i\le 4$. Letting $\mathcal{H}$ be the set of distinct components, we have by \eqref{7.7} and \eqref{7.34},
\[ \sum_{n=N+1}^{2N}\prod_{i=1}^4\Lambda_R(n, \mbf{H}_i) \ll_k
N \gs(\mathcal{H}) (\log R)^{\ell - |\mathcal{H}|}, \]
and by \eqref{2.13}
\[ \begin{split}\mathcal{T}(\ell_1,\ell_2,\ell_3,\ell_4) &\ll_k (\log R)^{4k-\ell}\sum_{j=1}^\ell \sum_{\substack{1\le h_1,\ldots ,h_j\le h\\ \text{distinct}}}N \gs(\mathcal{H}_j) (\log R)^{\ell - j} \\&  \ll_k N(\log R)^{4k} \sum_{j=1}^\ell \frac{h^j}{(\log R)^j}\\& 
\ll_k  N(\log N)^{4k},\end{split} \]
which proves \eqref{7.29}.

\end{document}